%
%
\font\bold=cmbx10 at 14pt

\font\tts=cmtt8
\centerline{\bold Homological Properties of the Homology Algebra of the Koszul }
\bigskip
\centerline{\bold Complex of a Local Ring. Examples and Questions.}

\medskip

\bigskip
\centerline{ Jan-Erik Roos}
\centerline{Department of Mathematics}
\centerline{Stockholm University}
\centerline{SE--106 91 Stockholm, SWEDEN}
\centerline{ e-mail: {\tt  jeroos@math.su.se}}
\bigskip
\rightline {\it Dedicated to J\"orgen Backelin at his $65^{th}$ birthday.}
\bigskip
\centerline{December 30, 2015}

\bigskip

\def\mysec#1{\bigskip\centerline{\bf #1}\nobreak\par}

\def\cite#1{~[{\bf #1}]}
\mysec{ Abstract.}
Let $R$ be a local commutative noetherian ring and $HKR$ the homology ring
of the corresponding Koszul complex. We study the homological properties of
$HKR$ in particular the so-called Avramov spectral sequence. When the 
embedding dimension of $R$ is four and when $R$ can be presented with
quadratic relations we have found 102 cases where this spectral sequence
degenerates and only three cases where it does not degenerate. We also
determine completely the Hilbert series of the bigraded Tor of these $HKR$
in tables A-D of section 5.  We also study
some higher embedding dimensions. Among the methods used are the programme
{\tt BERGMAN} by J\"orgen Backelin et al, the {\tt Macaulay2}-package {\tt DGAlgebras}
 by Frank Moore, combined with results by Govorov, Clas L\"ofwall, Victor Ufnarovski
 and others. 



\noindent  Hilbert series, Macaulay2, programme BERGMAN

\mysec{0. Introduction and Main Theorem}
Let $(R,m,k)$ be a local commutative noetherian ring with maximal ideal $m$ and residue field $k=R/m$.
Let ${\bf x}=(x_1,x_2,\ldots,x_n)$ be a minimal set of generators of the maximal ideal $m$. The Koszul complex
of $R$ is the 
exterior algebra $\oplus_{i=0}^n \Lambda^i R^n$ of a free $R$-module of rank n
equipped with the differential:
$$
d(T_{j_1}\wedge\ldots \wedge T_{j_i})=\sum_{l=1}^i(-1)^{l+1}x_{j_l}T_{j_1}\wedge \ldots \wedge {\widehat T_{j_l}}\wedge\ldots\wedge T_{j_i} 
$$
It will be denoted by $K({\bf x},R)$ or $KR$ since it is essentially independent of $\bf x$. It is a differential graded algebra
and its homology algebra $HKR$ is a skew-commutative graded algebra over $k$.
In the case that $R$ can be represented as a quotient of a regular local ring $(\tilde R,\tilde m)$ as
$R={\tilde R}/a$, where $a \subset {\tilde m}^2$ (passing to a completion of $R$ we can always assume that this is the case
for the problems we are studying), we have that $HKR$ is isomorphic to a {\rm Tor}-algebra:
$$
                     HKR \simeq {\rm Tor}_*^{\tilde R}({\tilde R}/a,k)
$$ 
This algebra $HKR$ has been studied in various special cases by many authors:
 it is an exterior algebra if and only if $R$ is a local complete
intersection; it is a Poincar\'e duality algebra it and only if $R$ is a Gorenstein local ring [A-G] and if $R$ is
a Golod ring then the square of the augmentation ideal  of $HKR$ is 0 (the converse is however not true,
even for rings with monomial relations [KAT]). But the general structure of $HKR$
(and in particular its homological properties) can be very complicated, even if $R$ is
a Koszul ring as we will see below. The aim of the present paper is to combine different methods towards studying $HKR$.
We will illustrate our methods on very precise examples, thereby discovering some new unexpected phenomena. We note that $HKR$ is graded, and one of our aims is to determine the double series
$$
\Phi_R(x,y)=\sum_{p\geq 0,q\geq 0}|{\rm Tor}^{HKR}_{p,q}(k,k)|x^py^q \leqno(1)
$$
(where $|V|$ denotes the dimension of a $k$-vector space 
$V$) for most
 (probably essentially all) quadratic rings $R$ of embedding dimension $4$:
Let $P_R(z)=\sum_{n\geq 0}|{\rm Tor}_n^R(k,k)|z^n$. We will see that with the exception of three explicit cases
we have $\Phi_R(z,z)=P_R(z)/(1+z)^4$. This last equality can be expressed by saying that the 
so-called Avramov spectral sequence {\it degenerates}.
Recall that this spectral sequence [AV1, formula (6.2.1)], is as follows for any local ring $R$
$$
E^2_n = \bigoplus_{p+q=n}E^2_{p,q}=\bigoplus_{p+q=n}{\rm Tor}^{HKR}_{p,q}(k,k) \Longrightarrow  {\rm Tor}^R(k,k) \otimes_{KR \otimes k} k 
$$
In particular there is a coefficientwise inequality $\ll$ of formal power series
$$
{P_R(z)\over (1+z)^m} \ll \Phi_R(z,z)
$$
where $m$ is the embedding dimension of $R$. 
 Of the three exceptional cases in embedding dimension 4, two are Koszul rings
with Hilbert series $1+4z+3z^2$ and $(1+3t-2t^3)/(1-t)$ respectively. In these two cases $HKR$ is
far from being a Koszul ring, indeed for the first case the bigraded Ext-algebra 
${\rm Ext}_{HKR}^*(k,k)$ has also generators in bidegrees $(3,8)$ and $(3,9)$.
Since there are other Koszul rings
 having the same series $1+4z+3z^2$
without this strange behaviour, we have here isolated new homological invariants even
for Koszul rings.
Let us recall that the Avramov spectral sequence always degenerates in embedding dimension $\leq 3$
and for higher embedding dimensions a consequence of it can be expressed by saing that $HKR$ and its so-called
Massey products determine $P_R(z)$. The converse is not true and the present paper shows that
by studying the homology of $HKR$ directly we sometimes obtain interesting invariants for $R$.
Our methods also show e.g. that for the ring $k[x,y,z,u]/(x^2,xy,yz,zu,u^2)$ 
the Avramov spectral sequence  {\it does} degenerate,
despite the fact that there is no differential algebra structure on the minimal resolution of the ideal
$(x^2,xy,yz,zu,u^2)$ over $k[x,y,z,u]$,
thereby answering a question by Avramov in [AV1] in the negative.
We will use a very useful package {\tt DGAlgebras} (written by Frank Moore [MOO]) which runs
under Macaulay2.
This package gives in particular an explicit presentation with generators and relations for $HKR$,
when you read in the ring $R$.
This presentation will then be analyzed, using results by Clas L\"ofwall [L1]
(the cube of the augmentation ideal of $HKR$ is often 0), J\"orgen Backelin et al [BA]
 (the programme {\tt BERGMAN}), Victor Ufnarovskij
(how to use {\tt BERGMAN} for determining Hilbert series of noncommutative graded algebras, where some of the generators
 have degrees $>1$) and others.
In the tables of Appendix 1 we have described $83$ different cases of the double series
$$
P_R(x,y)=\sum_{p\geq 0,q\geq 0}|{\rm Tor}^{R}_{p,q}(k,k)|x^py^q \leqno(2)
$$
 when $R$ is a quotient of
a polynomial ring $k[x,y,z,u]$ with an ideal generated by quadratic forms in $(x,y,z,u)$. (Note that the previous $P_R(z)=P_R(z,1)$.)
We have probably obtained all different such cases, but inside each case there are sometimes variations
due to the behavior of $HKR$, described in tables A,B,C,D of section 5 below. 
(this is like the periodic table of the elements,
where we also have isotopes). Examples of this are the cases {\bf 46}, {\bf 63} and {\bf 71}.  
More precisely here is our

MAIN THEOREM Let $R = k[x,y,z,u]/(f_1,f_2,\ldots, f_t)$, where $k$ is a field of characteristic $0$
 and the $f_i$:s are quadratic forms in the variables $x,y,z,u$ and $P_R(x,y)$ is given by (2) above.
Then there are $83$ different cases known for the $P_R(x,y)$ (they are given in Appendix 1).
Inside each of these cases $P_R(x,y)$ there are sometimes ``isotopes'' having different homological properties
of the homology algebra $HKR$ of the Koszul complex. In total there are $22$ such extra ``isotopes'' known, the
most extreme case of $P_R(x,y)$ is case 71 which is a Koszul algebra, but which has 4 extra isotopes.
With the exception of three of these $105=83+22$ cases we have that $\Phi_R(z,z)=P_R(z)/(1+z)^4$ (i.e. the so-called 
Avramov spectral sequence degenerates). The three exceptional cases being

a) Case 63ne (an isotope of Case 63) given by the ideal $(x^2,xy,xz+u^2,xu,y^2+z^2,zu)$ and where
       $$ (1+z)^4/P_{R_{63ne}}(z)-1/\Phi_{63ne}(z,z)=z^7(1+z)/(1-z+z^2)$$

b) Case 71v16 (an isotope of Case 71) given by the ideal $(xz+u^2,xy,xu,x^2,y^2+z^2,zu,yz)$ and where
 $$ (1+z)^4/P_{R_{71v16}}(z)-1/\Phi_{R_{71v16}}(z,z)=z^7(1+z)/(1-z+z^2)$$

c) Case 60va given by the ideal $(x^2+yz+u^2,xy,zu,z^2,xz+yu,xu)$ and where
 $$ (1+z)^4/P_{R_{60va}}(z)-1/\Phi_{R_{60va}}(z,z)=z^7(1+z)(1-3z^2-z^3)/(1-z+z^2)$$ 
   
Note that a) and b) are Koszul rings. We have also calculated the two-variable 
$\Phi_R(x,y)$ for the 105 cases.
The results are given in tables A-D section 5 and the methods used  for this are mentioned in section 5
In section 1 we start with case a) of the Main Theorem as a typical example
 of the calculation of $HKR$ and its homological properties.
In sections 2 and 3 we treat the cases c) and b) of out Main Theorem. Section 4 solves a
question of Avramov and in section 5
 we give the tables (tables A-D just mentioned)
 where hopefully all cases of $\Phi_{HKR}$ for the embedding dimension 4 cases (quadratic relations)
have been calculated. We assume  characteristic of $k$ to be $0$ throughout except in section 6 where 
we present conjectures and some results for the higher embedding dimension cases. 

We finally wish to thank Clas L\"ofwall, J\"orgen Backelin and Victor Ufnarovski for stimulating discussions.

\mysec{1. A Koszul ring in 4 variables with 6 relations which is exceptional}

Recall first that if $HKR$ is the homology of the Koszul complex of a ring $R$ of embedding dimension $m$
as in the introduction,
then there is a spectral sequence (the Avramov spectral sequence [AV1])
$$
E^2_n= \oplus_{p+q=n} {\rm Tor}_{p,q}^{HKR}(k,k)  => {\rm Tor}^R_n(k,k)/KR\otimes k
$$
leading to the inequality
$$
P_R(z)/(1+z)^m \leq \sum_{p\geq 0, q\geq 0,p+q=n}|{\rm Tor}_{p,q}^{HKR}(k,k)|z^n
$$
where we have equality if and only if the spectral sequence degenerates.

Let us now consider the ring 
$$
R=k[x,y,z,u]/(x^2,xy,zu,xu,y^2+z^2,xz+u^2)
$$. 
This is the case a) of the main theorem above.
The Koszul dual of this ring is $k<X,Y,Z,U>/(YU-UY,YZ-ZY,YY-ZZ,XZ-ZX-UU)$ and the Gr\"obner basis of this ideal is
{\it quadratic} if we use the ordering of the variables $Y,U,Z,X$. Thus $R$ is 
 a Koszul ring and its Hilbert series is $(1 + 3 t - 2 t^3)/(1 - t)$.
However, the homology of the Koszul complex of $R$, namely $HKR$ and its homology are rather complicated:
We now use the Macaulay2 package {\tt DGAlgebras} with the input file:

{\tts loadPackage(``DGAlgebras'')

R=QQ[x,y,z,u]/ideal(x*z+u\^{}2,x*y,x*u,x\^{}2,y\^{}2+z\^{}2,z*u)

res(ideal R); betti oo

res(coker vars R,LengthLimit => 7); betti oo

HKR=HH koszulComplexDGA(R)

generators HKR

for n from 1 to length(generators HKR) list degree X{\_}n

ideal HKR

I=ideal(vars HKR)

I\^{}2; trim(oo)

I\^{}3; trim(oo)

res(coker vars HKR,LengthLimit => 4);

betti oo}

From the output file of this we get first the important fact that the algebra $HKR$ has the cube of
its augmentation ideal $I=\overline{HKR}$ equal to 0. 
This first result is useful for  determining the homological properties of $HKR$, in particular
for determining the series
$$
\Phi_R(x,y)=\sum_{p,q \geq 0} |{\rm Tor}^{HKR}_{p,q}(k,k)|x^py^q
$$
Indeed, we can now use a general
 result by Clas L\"ofwall, which essentially says that many homological properties of a graded algebra,
whose augmentation ideal has cube 0 are determined by its Hilbert series and the subalgebra of the
whole {\rm Ext}-algebra generated
(under Yoneda product) of the ${\rm Ext}^1(k,k)$. More precisely, the theorem 1.3 on page
310 of [L1] (note that the citation page 19 in this Theorem should be page 309 {\it loc.cit.}). 
This theorem says that in our special case with the notations of {\it loc. cit.}  our
$\Phi_R(x,y)$ of (1) is equal to the  L\"ofwall $P_{HKR}(x,1,y,1)$) and the L\"ofwall formula gives
$$
\Phi_R(x,y)=xH_A(x,y)/(1+x-H_A(x,y)(1-H_{I/I^2}(y)x+H_{I^2}(y)x^2)) \leqno(3)
$$
where $H_A(x,y)$ is the Hilbert series in two variables of the algebra $A$ 
which is the subalgebra generated by ${\rm Ext}^1_{HKR}(k,k)$ of the Yoneda {\rm Ext}-algebra
 ${\rm Ext}^*_{HKR}(k,k)$.
This algebra is equal to the Koszul dual of $HKR$, which can be presented as 
$$
T((I/I^2)^*)/((im(I/I^2\otimes I/I^2 \rightarrow I^2)^*)
$$
where $T$ is the tensor algebra on $(I/I^2)^*$ which is the k-vector space dual of $I/I^2$ and $(im(I/I^2\otimes I/I^2 \rightarrow I^2)^*)$
is the twosided ideal in $T$ generated by the dual of the image of the natural multiplication map of $HKR$.
It is thus natural to write $HKR^!$ instead of $A$ and we will often write the preceding formula as
$$
1/{\Phi}_R(x,y)=(1+1/x)/HKR^!(x,y)-HKR(-x,y)/x \leqno(4)
$$
where $HKR(x,y)$ is the Hilbert series in two variables of $HKR$.
We now use this result and   
continue the analysis of the output file. 
The formula (3) shows that to get $\Phi_R(x,y)$ we need the Hilbert series $A(x,y)=HKR^!(x,y)$ of $A$
 and the two Hilbert series
$H_{I/I^2}(y)$ and $H_{I^2}(y)$, the most difficult part being the determination of  $A(x,y)$ and we start with that.
 We see that that we can determine a presentation of $HKR$ in our special
(and typical) case as having 
  $6$ generators $X_1,\ldots,X_6$ of degree $(1,2)$,
$6$ generators $X_7,\ldots,X_{12}$ of degree $(2,3)$ and one generator $X_{13}$
of degree $(3,4)$ (here the second element of these pairs of degrees comes from the grading of $R$ and will be ignored),
 satisfying the following {\it quadratic} relations:
$$
X_5X_6,X_4X_6,X_2X_6,X_1X_6,X_4X_5-X_3X_6,X_2X_5,X_1X_5,X_2X_4,X_1X_4,X_1X_3,X_1X_2,
$$
$$
X_6X_{12},X_5X_{12},X_2X_{12},X_1X_{12},X_6X_{11}-X_4X_{12},X_5X_{11}-X_3X_{12},X_1X_{11},X_6X_{10},
$$
$$
X_5X_{10},X_4X_{10},X_2X_{10},X_1X_{10},X_6X_{9},X_5X_{9},X_4X_{9},X_3X_{9}+X_4X_{12},X_2X_{9},
$$
$$
X_1X_{9},X_6X_{8},X_5X_{8}-X_4X_{12},X_4X_{8},X_3X_{8}+X_4X_{11},X_2X_{8},X_1X_{8},X_6X_{7},X_{5}X_{7},
$$
$$
X_{4}X_{7},X_3X_{7}+X_2X_{11},X_2X_{7},X_{1}X_{7},X_6X_{13},X_5X_{13},X_4X_{13},X_2X_{13},X_{1}X_{13},X_{12}X_{12},
$$
$$
X_{11}X_{12},X_{10}X_{12},X_9X_{12},X_8X_{12},X_{7}X_{12},X_{10}X_{11}+X_3X_{13},X_9X_{11},X_8X_{11},X_{7}X_{11},X_{10}X_{10},
$$
$$
X_{9}X_{10},X_{8}X_{10},X_7X_{10},X_9X_{9},X_{8}X_{9},X_{7}X_{9},X_8X_{8},X_7X_{8},X_7X_{7},X_{12}X_{13},X_{11}X_{13},
$$
$$
X_{10}X_{13},X_{9}X_{13},X_8X_{13},X_7X_{13}
$$

Now, using these relations we get that the Koszul dual $A$ of $HKR$ is the quotient of the free associative
algebra $k<Y_1,\ldots,Y_{13}>$ on the dual generators $Y_i$ of the $X_i$ with the ideal generated by
$$
[Y_2,Y_3],[Y_3,Y_4],[Y_3,Y_5],[Y_4,Y_5]+[Y_3,Y_6],[Y_3,Y_{10}],[Y_3,Y_{11}],Y_{11}Y_{11},[Y_3,Y_7]-[Y_2,Y_{11}],
$$
$$
[Y_5,Y_{11}]+[Y_3,Y_{12}],[Y_{10},Y_{11}]-[Y_3,Y_{13}],[Y_{3},Y_{8}]-[Y_4,Y_{11}],
$$
$$
[Y_4,Y_{12}]+[Y_6,Y_{11}]+[Y_5,Y_{8}]-[Y_3,Y_9]
$$
Here we have written e.g. $[Y_2,Y_3]=Y_2Y_3-Y_3Y_2$ and e.g. $[Y_3,Y_{13}]=Y_3Y_{13}-Y_{13}Y_3$ for the ``odd'' 
generators $Y_i,i=1\ldots 6$ and $Y_{13}$. Furthermore we have written e.g. $[Y_{10},Y_{11}]=Y_{10}Y_{11}+Y_{10}Y_{11}$
for the ``even'' generators and e.g. $[Y_3,Y_{10}]=Y_3Y_{10}+Y_{10}Y_3$ when the even and odd generators are
mixed.
Now to use the formula of L\"ofwall, we need to calculate the Hilbert series $A(x,y)$ of $A$ with 
{\tt BERGMAN} [BA]   using the 
presentation we have just detemined, where the bigrading of the preceding generators should
now be 1 for the first degree (corresponding to x) and for the second degree (corresponding to y) we should
have 1 for $Y_i (i=1,\ldots 6)$, 2 for $Y_i (i=7,\ldots 12)$ and 3 for $Y_{13}$.
  But {\tt BERGMAN} can for the moment only directly handle Hilbert series of algebras where the generators have only degree 1. 
 But if for the moment we are only interested in the generating series of the total terms
$$
\Phi_R(x,x)=\sum_{n=0}^{\infty}\sum_{p+q=n}|{\rm Tor}_{p,q}^{HKR}(k,k)|x^n
$$
which occur in the Avramov spectral sequence it is sufficient
to determine the Hilbert series $H_A(x,x)$ of A, i.e. to give the variables $Y_i (i=1,\ldots 6)$ the degree  2,
the variables $Y_i (i=7,\ldots 12)$ the degree 3 and the variable $Y_{13}$ the degree 4, keep the relations
as above and use the
following Lemma that Victor Ufnarovski has been kind to communicate to me:

LEMMA (Ufnarovski) Let $A$ be a graded algebra which is quotient of a free algebra with $n$ generators 
$Y_i$ of degrees $d_i \geq 1$ by an ideal $I=(f_1,\ldots,f_m)$ which is homogeneous with respect to the given grading. 
Let now $B$ be a new graded algebra obtained from $A$ as follows:
Introduce a new variable $t$ of degree 1 and new variables $Z_i$ of degrees 1, and let $B$ be the
quotient of the free associative algebra on the variables $t,Z_i$ with the ``same'' ideal as in $A$,
but where we have replaced each $Y_i$ by $t^{d_i-1}Z_i$. Then we have the following relation between
the Hilbert series of $A$ and $B$:
$$
1/A(z) = 1/B(z)+(n+1)z-\sum_{i=1}^n z^{d_i}
$$
i.e. we can reduce the calculation of $A(z)$ to the calculation of $B(z)$ which can be done directly in
{\tt BERGMAN} since all generators of $B$ have degree 1.

The proof of the Lemma will be given at the end of this section.

Remark.There is an analogous result in Ufnarovski's book [U], Lemma on page 141.   

We now continue with the proof of case a) of our main theorem:

We follow the recipe of  the Ufnarovski Lemma.
Recall that the variables $Y_i$ of degree (1,2) should now be considered to have degree 2 etc.

 We introduce a new variable $t$ of degree $1$ and replace the variables $Y_i$ above as follows:
$$
  Y_i=tZ_i,\,1\leq i \leq 6;\,Y_i=t^2Z_i,\,7\leq i \leq 12,\,Y_{13}=t^3Z_{13}
$$
where the new $Z_i$:s are given the degree $1$.
Now the algebra $A$ becomes the quotient of the free associative algebra in $14$ generators of degree $1$: 
$k<Z_1,\ldots,Z_6,t,Z_7\ldots,Z_{13}>$ with the ideal above where e.g. $[Y_2,Y_3]$ is replaced by
$tZ_2tZ_3-tZ_3tZ_2$ and  $[Y_{10},Y_{11}]-[Y_3,Y_{13}]$ is replaced by 
$$
t^2Z_{10}t^2Z_{11}+t^2Z_{11}t^2Z_{10}-tZ_3t^3Z_{13}+t^3Z_{13}tZ_3  
$$ 
etc. The new algebra is denoted by $B$ and the old one by $A$ and the Ufnarovski Lemma gives that we
following relation between their Hilbert series $H_A(z)$ and $H_B(z)$:
$$
1/H_A(z)=1/H_B(z)+14z-6z^2-6z^3-z^4 \leqno(5)
$$
It remains to calculate $H_B(z)$. For this we can now directly use {\tt BERGMAN} on the following input file
which we call {\tt inbcaseB} (we have written {\tts xi} instead of $Z_i$):

{\tts
(noncommify)

(algforminput)

vars x1,x2,x3,x4,x5,x6,t,x7,x8,x9,x10,x11,x12,x13;

t*x2*t*x3-t*x3*t*x2,t*x3*t*x4-t*x4*t*x3,t*x3*t*x5-t*x5*t*x3,

t*x4*t*x5-t*x5*t*x4+t*x3*t*x6-t*x6*t*x3,t*x3*t\^{}2*x10+t\^{}2*x10*t*x3,

t*x3*t\^{}2*x11+t\^{}2*x11*t*x3,

t\^{}2*x11*t\^{}2*x11,t*x3*t\^{}2*x7+t\^{}2*x7*t*x3-t*x2*t\^{}2*x11-t\^{}2*x11*t*x2,

t*x5*t\^{}2*x11+t\^{}2*x11*t*x5+t*x3*t\^{}2*x12+t\^{}2*x12*t*x3,

t\^{}2*x10*t\^{}2*x11+t\^{}2*x11*t\^{}2*x10-t*x3*t\^{}3*x13+t\^{}3*x13*t*x3,

t*x3*t\^{}2*x8+t\^{}2*x8*t*x3-t*x4*t\^{}2*x11-t\^{}2*x11*t*x4,

t*x4*t\^{}2*x12+t\^{}2*x12*t*x4+t*x6*t\^{}2*x11+t\^{}2*x11*t*x6+t*x5*t\^{}2*x8+

t\^{}2*x8*t*x5-t*x3*t\^{}2*x9-t\^{}2*x9*t*x3;
}

The following command in {\tt BERGMAN}:

{\tt (ncpbhgroebner ``inbcaseB'' ``t1'' ``t2'' ``t3'')}

\noindent gives after a few seconds {\tt t3}, i.e. the Hilbert series for $B$ in degrees $\leq 20$ and $\geq 2$, so that 
$$
1/H_B(z)=1-14z+4z^4+6z^5+2z^6-z^8-z^9+z^{11}+z^{12}-z^{14}-z^{15}+z^{17}+z^{18}-z^{20}+O(z^{21}) \leqno (6)
$$
leading to the guess (alternatively using the {\tt maple} command {\tt convert(``,ratpoly)}
on the previous formula (6))
that 
$$
1/H_B(z)={1-15z+15z^2-14z^3+4z^4+2z^5+4z^7+z^8 \over 1-z+z^2} \leqno(7)
$$
But this is not a proof. To get a proof we will analyze the file {\tt t1}, i.e. the groebner basis of the ideal in {\tt B}.
For the ordering of the variables we have given in {\tt inbcaseB} this groebner basis is infinite.
But, at my request some time ago J\"orgen Backelin wrote an addition {\tt permutebreak.sl} to the programme {\tt BERGMAN}
which checks the groebner bases up to a certain given degree of the different permutations of the variables of a given
input to the programme {\tt BERGMAN}. In the present case this seems rather discouraging in view of the fact that
$14! = 87178291200$ but rather immediately we see that switching the variables {\tt x2,x3} to {\tt x3,x2} in the {\tt vars}
line of {\tt inbcaseB} gives a finite groebner bases (in degrees $\leq 8$. Even better: after a short time with breaking
degree $12$, checking about $2700$ cases, we obtain that the ordering {\tt x1,x2,x4,x5,x6,t,x3,x7,x8,x9,x10,x11,x12,x13}
gives a finite groebner bases with $4$ elements of degree $2$, $6$ elements of degree $3$ and $2$ elements of degree $4$.
We therefore get $12$ leading monomials in the free algebra and the general structure of the Hilbert series follows from
Govorov ([GOV],Theorem 2) and the formula (7) follows. 
Now the formula (5) gives that
$$
1/A(z)={1-15z+15z^2-14z^3+4z^4+2z^5+4z^7+z^8 \over 1-z+z^2}+14z-6z^2-6z^3-z^4 =
$$
$$
= {1-z-5z^2+3z^4-3z^5-z^6+4z^7+z^8 \over  1-z+z^2}
$$
and now to use the L\"ofwall formula we only need to observe that
$I^2$ is equal to
$$
X_3X_6,X_3X_5,X_3X_4,X_2X_3,X_4X_{12},X_3X_{12},X_4X_{11},X_3X_{11},X_2X_{11},X_3X_{10},X_3X_{13},X_{11}X_{11}
$$
and $I^3=0$ so that 
$$
HKR(x,y)=1+6xy+6xy^2+xy^3+4x^2y^2+6x^2y^3+2x^2y^4
$$
and 
$$
1/\Phi_R(z,z)=(1+1/z)/A(z)-(1-6z^2-6z^3+3z^4+6z^5+2z^6)/z=
$$
$$
={(1+z)(1-2z-3z^2+3z^3-3z^5+2z^6+z^7)\over 1-z+z^2}
$$

Now, since $R$ is a Koszul ring with Hilbert series $P_R(z)={1+z \over 1-3z+2z^3}$ given above, it follows in particular that {\it if} the
Avramov spectral sequence degenerated, we would have had $P_R(z)/(1+z)^4 = \Phi_R(z)$. But
$$
(1+z)^4/P_R(z)-1/\Phi_R(z) = {z^7(1+z)\over 1-z+z^2}
$$
shows that this is not the case (there is a non-zero Massey product in $KR$ of degree $7$).

Remark: In the last section we will see how to calculate the series of $\Phi_R(x,y)$ using the extra
degree in $HKR$ coming from the Koszul homology. There is also a third degree, coming from the
grading of $R$, that we will not use.

PROOF OF THE UFNAROVSKI LEMMA (UFNAROVSKI):
We need the theory of $n$-chains [U] of an algebra such as $A=k<Y_1,Y_2,\ldots>/(f_1,\ldots f_m)$.
We use deglex order.
It is known (cf. [U],p.51) that
$$
1/A(z) = 1-H_X(z)+H_{C_1}(z)-H_{C_2}(z)+H_{C_3}(z)- \ldots
$$
where $X={Y_1,Y_2,\ldots}$ are the generators of $A$, and 
$H_X(z)=\sum z^{d_i}$ is the Hilbert series (i.e. the generating series) of X, and where
$C_1$=the leading terms of the relations $f_i$ (we use the deglex order)
and $H_{C_1}(z)$ is the generating series of $C_1$, similarly for the $n$-chains $C_n$ for $n\geq 2$.
We now introduce an auxiliary algebra $D$ which I write 
$$
D=<Y_1,Y_2,\ldots,t,Z_1,Z_2,\ldots>/(f_1,\ldots,f_m,Y_i-t^{d_i-1}Z_i \ldots>
$$
Now if the order is deglex such that
$$
Y_1 > Y_2 > \ldots >t > Z_1 > Z_2 >\ldots Z_n
$$
For this algebra $D$ we have the new $X'=X \bigcup \{t,Z_1,\ldots,Z_n\}$ and the new 
$$
C_1'=C_1 \bigcup \{t^{d_i-1} Z_i\}
$$
 but the higher $C_i$ for $i\geq 2$ are the same.
It follows that 
$$
H_D(z)=H_A(z)-(n+1)z+\sum_{i=1}^n z^{d_i}
$$
On the other hand, consider the order
$$
Y_1>Y_2> \ldots Z_1 >Z_2 > \ldots >t
$$
We have that $Y_i$ can be replaced by $t^{d_i-1}Z_i$ 
and we get the algebra 
$$
B = k<Z_1,Z_2,\ldots Z_n,t>/(g_1,g_2,\ldots,g_m)
$$
where the relations $g_i$ are  obtained from the $f_i$ by replacing all $Y_i$ by $t^{d_i-1}Z_i$
Clearly $H_B(z)=H_D(z)$ and the lemma is proved.

\mysec{2. One more exceptional case.}
This  case will be treated as in section 1 so we can be more brief here.
The is the case c) in the main Theorem of section 0 i.e. the ring $R_{60va}=k[x,y,z,u]/(x^2+yz+u^2,xz+yu,xu,x^2+xy,z^2,xu)$.
 The resolution of the ideal of $R_{60va}$ over $k[x,y,z,u]$ has betti numbers:
$$\left ( \matrix{
 1 & .&.&.&.\cr
 .&6&5
&.&. \cr
 .&.&6 &9&3 \cr} \right )
$$
and an applications of the programme {\tt DGAlgebras} gives again that the algebra $HKR_{60va}$ has 
the cube of the augmentation ideal $\overline{HKR_{60va}}$ equal to 0. Furthermore we have
6 generators $X_1,\ldots,X_6$ of degree $(1,2)$ and 5 generators 
$X_7,\ldots,X_{11}$ of degree $(2,3)$ satisfying the following
{\it quadratic} relations:
$$
X_5X_6,X_3X_6+X_4X_6,X_1X_6,X_4X_5-X_2X_6,X_3X_5+X_2X_6,X_2X_5,X_1X_5,X_1X_4,X_1X_2,
$$
$$
X_6X_{11},X_5X_{11},X_3X_{11}+X_4X_{11},X_1X_{11},X_6X_{10},X_5X_{10},X_4X_{10}-X_2X_{11},X_2X_{10},
$$
$$
X_1X_{10},X_6X_{9}-X_4X_{11},X_5X_{9}-X_2X_{11},X_1X_{9},X_6X_{8},X_5X_{8},X_3X_{8}+X_4X_{8},
$$
$$
X_1X_{8}+X_3X_{10}+X_2X_{11},X_6X_{7}-X_2X_{11},X_5X_{7},X_4X_{7}+X_2X_{9},X_2X_{7},X_1X_{7},X_{11}X_{11},
$$
$$
X_{10}X_{11},X_9X_{11},X_8X_{11},X_7X_{11},X_{10}X_{10},X_9X_{10},X_8X_{10},X_7X_{10},X_7X_{9},X_{8}X_{8},X_{7}X_{7}
$$
Now, using these relations we get that the Koszul dual of $HKR$ is the quotient of the free algebra
$k<Y_1,\ldots,Y_{11}>$ on the dual generators $Y_i$ of the $X_i$ with the ideal generated by
$$
[Y_1,Y_3],[Y_2,Y_3],[Y_3,Y_4],[Y_2,Y_4],-[Y_3,Y_6]+[Y_4,Y_6],[Y_4,Y_5]+[Y_2,Y_6]-[Y_3,Y_5],[Y_2,Y_8],
$$
$$
[Y_3,Y_7],[Y_3,Y_9],[Y_4,Y_9],[Y_7,Y_8],[Y_8,Y_9],Y_9Y_9,-[Y_3,Y_{11}]+[Y_4,Y_{11}]+[Y_6,Y_9],
$$
$$
[Y_4,Y_{10}]+[Y_2,Y_{11}]+[Y_5,Y_9]-[Y_1,Y_8]+[Y_6,Y_7],-[Y_3,Y_8]+[Y_4,Y_{8}],-[Y_1,Y_{8}]+[Y_3,Y_{10}],
$$
$$
-[Y_4,Y_{7}]+[Y_2,Y_{9}] \leqno(9)
$$
Here we have again written e.g. $[Y_1,Y_3]=Y_1Y_3-Y_3Y_1$ for the ``odd'' generators $Y_i,i=1\ldots 6$,
and e.g. $[Y_8,Y_9]=Y_8Y_9+Y_9Y_8$ for the ``even'' generators $Y_i,i=7\ldots 11$, and e.g.
$[Y_2,Y_8]=Y_2Y_9+Y_8Y_2$ when even and odd generators are mixed.
We want to calculate the Hilbert series of this quotient, when the generators $Y_i$ are given
the degree 2 for $i=1,\ldots,6$ and degree 3 for $i=7,\ldots,11$. For this we use again the
programme {\tt BERGMAN} and the idea of Victor Ufnarovski to handle the gradings, i.e. we introduce
 a new variable $t$ of degree 1 and replace the variables $Y_i$ in (9)
$$
Y_i=tZ_i , 1 \leq i \leq 6; Y_i=t^2Z_i, 7 \leq i \leq 11
$$
where the $Z_i$:s also have degree 1. Now the algebra $A$ becomes the quotient of the free algebra in 
12 generators of degree 1: $k<Z_1,\ldots,Z_6,t,Z_7,\ldots,Z_{11}>$
with the ideal (9) above where e.g. $[Y_1,Y_3]$ is replaced by $tZ_1tZ_3-tZ_3tZ_1$ and e.g. $Y_9Y_9$ is replaced
by $t^2Z_9t^2Z_9$ and $-[Y_4,Y_7]+[Y_2,Y_9]$ is replaced by 
$$
-tZ_4t^2Z_7-t^2Z_7tZ_4+tZ_2t^2Y_9+t^2Y_9tZ^2
$$
etc.  If the new algebra is
denoted by $B$ and the old one by $A$, then by the theory of Ufnarovski we have the following relation
between their Hilbert series:
$$
1/H_A(z)=1/H_B(z)+12z-6z^2-5z^3
$$
It remains to calculate $H_B(z)$.

After these preparations we are now ready again to use the programme {\tt BERGMAN} and the command in {\tt BERGMAN}

 {\tt(ncpbhgroebner} $\,''${\tt inbcaseB}$''\quad''${\tt t1}$''\quad''${\tt t2}$''\quad''${\tt t3}$''\,${\tt )}

After a few seconds we get  {\tt t3} (i.e. the Hilbert series for $B$) in degrees $\leq 20$ and $\geq 2$ and
we obtain with {\tt maple}
$$
S:= {\tt series}(1/(1+12z+{\tt t3},z,21) = 1-12z+6z^4+9z^5+2z^6-3z^7-4z^8-2z^9+z^{10}
$$
$$
+3z^{11}+2z^{12}-z^{13}-3z^{14}-2z^{15}+z^{16}+3z^{17}+2z^{18}-z^{19}-3z^{20}+O(z^{21})
$$
and the maple command {\tt convert(S, ratpoly)} on the previous series gives the rational function:
$$
(1-z)(z^9  + 2 z^8  + z^7  - 3 z^6  - 2 z^5  - 5 z^4  - 11 z^3  + z^2  - 12 z + 1)\over(1-z+z^2) \leqno(9)
$$
 Here is again how one should proceed to get a proof of this formula.
If we can find a permutation of the 12 variables in {\tt inbcaseB} so that the groebner basis {\tt t1} is
finite, then the Hilbert series of $B$ is equal to the Hilbert series of the free algebra
on the 12 variables divided by the ideal generated by the monomial leading terms of this
finite groebner basis, and this last Hilbert series is rational of a special form according to
Govorov. Again we use J\"orgen Backelin's addition {\tt permutebreak.sl} to the programme {\tt BERGMAN}
and we obtain that the order of the variables 
$x3,t,x1,x2,x5,x6,x4,y8,y7,y9,y10,y11$
gives a finite groebner basis with 6 elements of degree 4, 9 elements of
degree 5 and 3 elements of degree 6. We therefore get 18 leading monomials
in the free algebra and the general structure of the Hilbert series 
follows from Govorov ([GOV], Theorem 2) and the formula (9) follows.
We have now all we need to prove the following

THEOREM  Let $R_{60va}=k[x,y,z,u]/(x^2+yz+u^2,xz+yu,zu,x^2+xy,z^2,xu)$. Then

i)
$$
\Phi_{R60va}(z,z)=
$$
$$\sum_{n=0}^{\infty}\sum_{p+q=n}|{\rm Tor}^{HKR}_{p,q}(k,k)|z^n =
{(1-z+z^2)\over(1+z)(z^9 - z^8 + z^7 + z^6 - 4z^5 + z^4 + 4z^3 - 3z^2 - 2z + 1)}
$$ 
ii) $P_R(z)=(1+z)/(-z^6 - z^5 + 3z^3 - 3z + 1)$

iii)$ 1/\Phi_{R60va}(z,z)-(1+z)^4/P_{R60va}(z) = z^7(1+z)(1-3z^2-z^3)/(1-z-z^2)$

It follows that the Avramov spectral sequence does not degenerate (the two series differ starting in degree 7).
\medskip
PROOF: For the algebra $HKR$ we have that the cube of the augmentation ideal $\overline{HKR}$ is zero.
Furthermore the rational function $1/HKR^!(z)$ is given by the formula (9) above and the two-variable form
of the Hilbert series of $HKR$ is given by the formula
$$
HKR(x,y)=1+6xy+5xy^2+6x^2y^2+9x^2y^3+3x^2y^4
$$
This gives the one-variable Hilbert series  $H(-z,z)$
and the formula
$$
1/P_{HKR}(z) = (1+1/z)/HKR^!(z)-H(-z,z)/z
$$
gives the result i).
\medskip
To prove ii) we first observe that the 
$R^!=k<X,Y,Z,U>/(Y^2,X^2-XY-YX-YZ-ZY,XZ+ZX-YU-UY,YZ+ZY-U^2)$,
where the $X,Y,Z,U$ are the dual variables to $x,y,z,u$.
There are only 24 permutations of the variables $X,Y,Z,U$ and one
sees that the order $X,Y,U,Z$ gives in BERGMAN a finite Gr\"obner basis,
and that $1/R^!(z)=(1-z)^2(1-2z-z^2-z^3)$. On the other hand the Hilbert
series $R(z)=(1-3z+3z^3)/(1-z)=1+4z+4z^2+z^3+z^4+z^5+\ldots$
and according to Theorem B.9 by Clas L\"ofwall (cf. page 310 of [R3] the ring $R$
satisfies a condition ${\cal M}_3$ so that we still have the formula
$1/P_R(z)=(1+1/z)/R^!(z)-R(-z)/z$ and this gives ii) and the assertion iii) follows.

\mysec{3. One final exceptional case for a local ring that is artinian and a Koszul algebra.}
In section 4 we will treat a  Koszul local ring (case 46) whose Koszul algebra has properties
that show that the Avramov spectral sequence degenerates (solving a problem of Avramov)
In section 2 we found a quadratic local ring (case 60va) whose Avramov spectral sequence does not degenerate.
But case 60va is not a Koszul algebra.
Now we pass to a case of a Koszul local ring (case 71) which has further rather unexpected
properties and which is case b) of the Main Theorem.
\medskip
THEOREM.- There are at least five local commutative rings having Hilbert series $(1+z)(1+3z)$
which are  Koszul algebras and which have different homological properties.
Furthermore for one of these five cases the Avramov spectral sequence does not degenerate.
\medskip
PROOF: For the ring of case 71 in the Appendix 1 we have the presentation
$$R_{71}=k[x,y,z,u]/(x^2,y^2,z^2,u^2,xy,zu,yz+xu)$$
But we also have three other rings which are Koszul algebras and have the same Hilbert series:
$$ 
{\scriptstyle R_{71v4}={k[x,y,z,u]\over(x^2+u^2,xy,xu,y^2,yz,z^2,zu)}} \, {\scriptstyle R_{71v7}={k[x,y,z,u]\over(x^2,y^2,z^2,xz+u^2,xu,yz,zu)}}\, {\scriptstyle R_{71v5}={k[x,y,z,u]\over(x^2+xy,x^2+yz,xy+y^2,z^2,xu,zu,xz+u^2)}}
$$
and but the betti numbers of the resolution of the preceding {\it four} ideals are:
$$
\left ( \matrix{
 1 & .&.&.&.\cr
 .&7&8
&2&. \cr
 .&.&5 &8&3 \cr} \right ) \quad
\left ( \matrix{
 1 & .&.&.&.\cr
 .&7&8
&3&. \cr
 .&.&6 &8&3 \cr} \right ) \quad
\left ( \matrix{
 1 & .&.&.&.\cr
 .&7&8
&1&. \cr
 .&.&4 &8&3 \cr} \right ) \quad
\left ( \matrix{
 1 & .&.&.&.\cr
 .&7&8
&.&. \cr
 .&.&3 &8&3 \cr} \right ) 
$$
so that the rings can not be isomorphic. But their $HKR$ are rather easy, in one of the cases $HKR$ has
only monomial relations, But strangely enough the case of a $R_{71}$ with only monomial relations
gives a more complicated $HKR$.

But there remains one case $R_{71v16}$ which is {\it very} different from the preceding and leads
to a new phenomenon. We now treat it in detail:
The ring $$
R_{71v16}=k[x,y,z,u]/(xz+u^2,xy,xu,x^2,y^2+z^2,zu,yz)
$$
is also a Koszul algebra and it has the same matrix of Betti numbers as $R_{71}$
 but the algebraic behavior of its $HKR$ is very different.
It still has the cube of the augmentation ideal of $HKR$ equal to 0, it has 
seven $X_1,\ldots,X_7$ generators of degree $(1,2)$, eight generators $X_8,\ldots,X_{15}$ of degree $(2,3)$
 and two generators $X_{16},X_{17}$ of degree $(3,4)$ 
and an applications of the programme {\tt DGAlgebras} gives that these generators satisfy
$128$ {\it quadratic relations} of which however only the following $16$ are non-monomial relations:
$$
X_5X_6-X_3X_7,X_2X_3-X_4X_5,X_7X_{13}-X_5X_{15},X_6X_{13}-X_3X_{15},X_3X_{12}-X_5X_{14},
$$
$$
X_3X_{11}+X_5X_{15},X_6X_{10}-X_5X_{15},X_4X_{10}-X_2X_{13},X_3X_{10}+X_5X_{13},X_5X_9+X_2X_{13},
$$
$$
X_3X_{9}+X_4X_{13},X_3X_{8}+X_2X_{13},X_3X_{16}-X_5X_{17},X_{13}X_{14}+X_3X_{17},X_{10}X_{14}+X_5X_{17},
$$
$$
X_{12}X_{13}+X_5X_{17}
$$
and since the following 6 quadratic monomials do {\it not} occur among the 128 relations:
$$
X_3X_4,\,X_3X_5,\,X_3X_6,\,X_3X_{13},\,X_3X_{14},\,X_{13}X_{13}
$$
it follows that the Koszul dual $HKR^!$ is equal to the quotient of the free
associative algebra on the 17 variables variables $Y_i,1\leq i \leq 17$ 
(which are dual to the $X_i$) with respect
to the ideal generated by the 16 quadratic relations:
$$
[Y_3,Y_4],[Y_3,Y_5],[Y_3,Y_6],[Y_3,Y_{13}],[Y_3,Y_{14}],Y_{13}Y_{13},[Y_2,Y_3]+[Y_4,Y_5],[Y_5,Y_6]+[Y_3,Y_7]
$$
$$
[Y_3,Y_{12}]+[Y_5,Y_{14}],[Y_3,Y_{10}]-[Y_5,Y_{13}],[Y_3,Y_{15}]+[Y_6,Y_{13}],[Y_3,Y_9]-[Y_4,Y_{13}],
$$
$$
[Y_3,Y_{17}]-[Y_{13},Y_{14}],[Y_2,Y_{13}]-[Y_3,Y_8]+[Y_4,Y_{10}],[Y_5,Y_{15}]+[Y_7,Y_{13}]-[Y_6,Y_{10}]-[Y_3,Y_{11}],
$$
$$
[Y_5,Y_{17}]-[Y_{12},Y_{13}]-[Y_{10},Y_{14}]+[Y_3,Y_{16}]
$$
Here $[,]$ means as before the graded commutator so that $[Y_i,Y_j]=Y_iY_j-Y_jY_i$ if $1\leq i \leq 7$ or $i=16$ or $i=17$
and $1\leq j \leq 7$ or $j=16$ or $j=17$. Furthermore we have $[Y_s,Y_t]=Y_sY_t+Y_tY_s$ otherwise.
We now need to calculate the Hilbert series of $HKR^!$ where we have given the variables $Y_1,\ldots,Y_7$ the degrees 2,
the variables $Y_8,\ldots,Y_{15}$ the degrees 3 and the variables $Y_{16},Y_{17}$ the degrees 4. We use the Ufnarovski
trick and introduce a new variable $T$ of degree 1 and replace everywhere in $HKR^!$ the variables $Y_i$ according to
the formulae $Y_i=TZ_i$ for $1\leq i \leq 7$ , $Y_i=T^2Z_i$ for $8\leq i \leq 15$ and $Y_i=T^3Z_i$ for i=16,17. 
The modified algebra $A$ has generators of degree 1 and can be put into {\tt BERGMAN} and its Hilbert
series $A(z)$ is related to $HKR^!(z)$ by the Ufnarovski formula:
$$
{1\over HKR^!(z)} = {1\over A(z)}+18z-7z^2-8z^3-2z^4
$$
The calculations in {\tt BERGMAN} gives after a few minutes $A(z)$ up to degree 22, giving  that
$$
1/A(z)=1-18z+5z^4+8z^5+3z^6-z^8-z^8+z^{11}+z^{12}-z^{14}-z^{15}+z^{17}+z^{18}-z^{20}-z^{21}+O(z^{23})
$$
indicating (using the {\tt maple} command {\tt convert(-,ratpoly)} as earlier) that
$$
{1\over A(z)}={1-19z+19z^2-18z^3+5z^4+3z^5+5z^7+2z^8 \over 1-z+z^2}
$$
To prove this we use again the {\tt permutebergman} addition by
J\"orgen Backelin to {\tt BERGMAN} and we try to find a permutation of
the variables $T,Z_i,i=1\ldots 17$ that gives a finite Gr\"obner basis
when we calculate the Hilbert series of the algebra $A$. Since the number of
permutations of 18 variables is 18!=6402373705728000 this seems to be a
hopeless task, but after letting the programme run through only the 24000 first permutations
(this took about an hour) we finally got the permutation
$$
Z1,Z3,Z4,Z6,Z7,T,Z2,Z5,Z8,Z9,Z10,Z11,Z12,Z13,Z14,Z15,Z16,Z17
$$
for which the Gr\"obner basis is finite: it has 5 elements of degree 4, 8 elements of degree
5 and 3 elements of degree 6. Therefore, according to the result of Govorov
the Hilbert series $A(z)$ is rational as indicated.
We have now all we need to prove the following

THEOREM  Let $R=k[x,y,z,u]/(xz+u^2,xy,xu,x^2,y^2+z^2,zu,yz)$. Then
$R$ is a Koszul algebra with Hilbert series $R(z)=1+4z+3z^2$ so that
$P_R(z)=1/(1-4z+3z^2)$. But $HKR$ is far from being a Koszul algebra.
Indeed:

a)
$$
\Phi_R(z,z) =\sum_{n=0}^{\infty}\sum_{p+q=n}|{\rm Tor}^{HKR}_{p,q}(k,k)|z^n =
$$
$$
{(1-z+z^2)\over(1+z)(2z^7 + 2z^6 - 4z^5 + z^4 + 3z^3 - 4z^2 - 2z + 1)}
$$
b) Furthermore $P_R(z)=1/(1-4z+3z^2)$ so that
$$
(1+z)^4/P_R(z) -1/\Phi_R(z,z)=z^7(1+z)/(1-z+z^2)
$$
It follows again that the Avramov spectral sequence does not degenerate. 

PROOF: We have just determined $1/A(z)$ and it follows that 
$$
{1\over HKR^!(z,z)}={1-z-6z^2-z^3+4z^4-3z^5-2z^6+5z^7+2z^8 \over 1-z+z^2}
$$
But the two-variable version of the Hilbert series of $HKR$ is
$$
HKR(x,y)=1+7xy+8xy^2+2x*y^3+5x^2y^2+8x^2y^3+3x^2y^4
$$
This gives again a one-variable
$HKR(-z,z)=1-7z^2-8z^3+3z^4+8z^5+3z^6$
The formula by Clas L\"ofwall
$1/P_{HKR}(z,z)=(1+1/z)/HKR^!(z,z)-HKR(-z,z)/z$
now gives
$$
1/P_{HKR}(z,z)={(1+z)(1-2z-4z^2+3z^3+z^4-4z^5+2z^6+2z^7) \over 1-z+z^2}
$$
and the Theorem is proved.

\mysec{4. The other embedding dimension four cases.}

We have just treated the three interesting cases of the Main Theorem in section 0.
For the other cases we have that the Avramov spectral sequence {\it does} degenerate and
they are treated in an analogous way. Here we will only treat case 46 in detail
because it solves a problem that was left open in Avramov [AV1].
The case 46 in the ring $k[x,y,z,u]/(x^2,xy,yz,zu,u^2)$.

Application as before of the programme {\tt DGAlgebras} gives that the algebra $HKR$ has the cube of the 
augmentation ideal equal to 0 and that $HKR$ has 5 generators $X_1,\ldots,X_5$ of degree (1,2) and
4 generators $X_6,X_7,X_8,X_9$ of degree (2,3). These generators have only quadratic relations and they are:
$$
X_4X_5,X_3X_5,X_2X_5,X_3X_4,X_2X_3,X_1X_3,X_1X_2,X_5X_9,X_4X_9,X_3X_9,X_2X_9,
$$
$$
X_5X_8,X_4X_8,X_3X_8,X_5X_7,X_4X_7-X_2X_8,X_3X_7,X_2X_7,X_1X_7,X_5X_6-X_1X_9,
$$
$$
X_3X_6,X_2X_6,X_1X_6,X_9X_9,X_8X_9,X_7X_9,X_6X_9,X_8X_8,X_7X_8,X_7X_7,X_6X_7,X_6X_6
$$
so that the Koszul dual $HKR^!$ is the quotient of the free algebra
$k<Y_1,\ldots,Y_9>$ in the dual generators $Y_i$ of the $X_i$ with the ideal generated by
$$
[Y_1,Y_4],[Y_1,Y_5],[Y_2,Y_4],[Y_1,Y_8],[Y_4,Y_6],[Y_6,Y_8],[Y_4,Y_7]+[Y_2,Y_8],[Y_5,Y_6]+[Y_1,Y_9]
$$ 
Here again we have written e.g. $[Y_1,Y_4]=Y_1Y_4-Y_4Y_1$ for the ``odd'' generators and  e.g. $[Y_6,Y_8]=Y_2Y_8+Y_8Y_2$
for the ``even'' generators and e.g. $[Y_4,Y_7]=Y_4Y_7+Y_7Y_4$ when the even and odd generators are mixed.
We use the previous methods to calculate the Hilbert series of our quotient when the generators $Y_1,\ldots,Y_5$
are given the degrees 2 and the generators $Y_6,\ldots,Y_9$ are given the degree 3, i.e. we introduce a new variable
$t$ of degree 1 and replace the variables $Y_i$ by $tZ_i$ for $1\leq i\leq 5$ and $t^2Z_i$ for $6\leq i\leq 9$.
where the $Z_i$ now have degree 1. Our algebra $HKR^!$ is now replaced by the quotient of 
$k<Z_1,\ldots,Z_5,t,Z_6,\ldots,Z_9>$ by the ideal of $HKR^!$ where e.g. $Y_1Y_4-Y_4Y_1$ is replaced by
$tZ_1tZ_4-tZ_4tZ_1$ etc. If this new algebra is denoted by $B$ then again by the theory of Ufnarovskij
we have in this special case the following relation between the their Hilbert series $H_A(z)$ and $H_B(z)$
$$
1/H_A(z) = 1/H_B(z)+10z-5z^2-4z^3
$$
and $H_B(z)$ is calculated by {\tt BERGMAN} as before. We get:
$$
1/H_B(z)=1-10z+3z^4+4z^5+z^6
$$
since the permutation $Z_1,Z_2,Z_3,t,Z_4,\ldots,Z_9$   of the variables in $B$
gives a finite groebner basis (3 elements of degree 4, 4 elements of degree 5 and one element of degree 6).

We are now ready to calculate $\Phi(x,x)$ of our $HKR$. We have the formula of L\"ofwall
$$
1/\Phi_R(z)=(1+1/z)/A(z)-HKR(-z,z)/z
$$
where we have just found that 
$$
1/A(z)=(1-10z+3z^4+4z^5+z^6)+10z-5z^2-4z^3 =1-5z^2-4z^3+3z^4+4z^5+z^6
$$
For the Hilbert series in two variables of $HKR$ we observe that since
$$
{\overline {HKR}}^2=X_1X_5,X_2X_4,X_1X_4,X_1X_9,X_2X_8,X_1X_8,X_4X_6,X_6X_8
$$ 
that
$$
HKR(x,y)=1+5xy+4xy^2+3x^2y^2+4x^2y^3+x^2y^4
$$
so that $HKR(-z,z)=1-5z^2-4z^3+3z^4+4z^5+z^6$
so that $1/\Phi_R(z)=1-5z^2-4z^3+3z^4+4z^5+z^6$ too. But since $R$ is a Koszul ring (quadratic monomial relations)
with Hilbert series $H_R(z)=(1+3z+z^2-z^3)/(1-z)$ we have that (Fr\"oberg) $P_R(z)=1/H_R(-z)$ 
we have $(1+z)^4/P_R(z)=1-5z^2-4z^3+3z^4+4z^5+z^6$ too and therefore the Avramov spectral sequence degenerates.
 
\mysec{5. More precise $\Phi_R(x,y)$}

In the preceding sections we have shown how to calculate $\Phi_R(x,x)$
for some examples $R$.
In  this section we will given an indication about how to generalize this to the calculation
of the double series $\Phi_R(x,y)$. The complete calculation results for
the quadratic embedding dimension $\leq 4$ case are given in Tables A-D below.
Let us illustrate these more precise calculations with the simple example of case 46 which we have just treated for
the calculation of $\Phi_R(x,x)$.
 Note first that the output file for the case 46 gives for the dimensions
of the ${\rm Tor}_{p,q}^{HKR}(k,k)$ the integers:
$$
\halign{\tt#\hfill\quad&&\hfill#\quad\cr
0: & 1& 5& 22&  95&   409&  1760&  7573&  32585\cr
1: & .& 4& 36& 236&  1364&  7368& 38152& 191908\cr
2: & .& .& 15& 198&  1723& 12438& 80628& 487202\cr
3: & .& .&  .&  56&   976& 10576& 91448& 690904\cr
4: & .& .&  .&   .&   209&  4527& 58685& 590894\cr
5: & .& .&  .&   .&     .&   780& 20196& 304696\cr
6: & .& .&  .&   .&     .&     .&  2911&  87692\cr
7: & .& .&  .&   .&     .&     .&     .&  10864\cr
}
$$
showing that $\Phi_R(x,y)$ should start as
$$
\Phi(x,y)=1+xy(5+4y)+x^2y^2(22+36y+15y^2)+x^3y^3(95+236y+198y^2+56y^3)+\ldots
$$
But note that the calculation of $\Phi_R(z,z)$ in the preceding section gives
$$
\Phi_R(z,z)=1+5z^2+4*z^3+22z^4+36z^5+110z^6+236z^7+607z^8+1420z^9+3483z^{10}+\ldots
$$
and as you see this is not sufficient to get the preceding columns e.g. to find the decomposition 
$110=15+95=|{\rm Tor}_{2,4}^{HKR}(k,k)|+|{\rm Tor}_{3,3}^{HKR}(k,k)|$ etc.
But now we use the L\"ofwall formula (4) where we replace $x$ by $z^u$ for any integer $u \geq 1$ and $y$ by $z$ :
$$
1/\Phi_R(z^u,z)=(1+1/z^u)/H_A(z^u,z)-(1-H_{I/I^2}(z)z^u+H_{I^2}(z)z^{2u})z^u \leqno(10)
$$
gives that $H_A(z^u,z)$ is the Hilbert series of a new variant of the algebra $A$ where all previous generators of 
degree $s$ have been replaced by generators of degree $s+u-1$.
Now the Hilbertseries of $A$ with generators of these degrees can still be calculated
with {\tt BERGMAN}. We illustrate this with the case 46 just treated when $u=31$
The new input file for {\tt BERGMAN} will therefore be:

{\tts
(noncommify)

(algforminput)

vars x1,x2,x3,x4,x5,t,x6,x7,x8,x9;

t\^{}31*x1*t\^{}31*x4-t\^{}31*x4*t\^{}31*x1,t\^{}31*x2*t\^{}31*x4-t\^{}31*x4*t\^{}31*x2,

t\^{}31*x1*t\^{}31*x5-t\^{}31*x5*t\^{}31*x1,t\^{}31*x1*t\^{}32*x8+t\^{}32*x8*t\^{}31*x1,

t\^{}31*x4*t\^{}32*x6+t\^{}32*x6*t\^{}31*x4,t\^{}32*x6*t\^{}32*x8+t\^{}32*x8*t\^{}32*x6,

t\^{}31*x4*t\^{}32*x7+t\^{}32*x7*t\^{}31*x4+t\^{}31*x2*t\^{}32*x8+t\^{}32*x8*t\^{}31*x2,

t\^{}31*x5*t\^{}32*x6+t\^{}32*x6*t\^{}31*x5+t\^{}31*x1*t\^{}32*x9+t\^{}32*x9*t\^{}31*x1;
}

Note that here the generators still all have degree 1 and {\tt BERGMAN} now gives for this infile
the inverse of its Hilbert series:
$$
1-10*z+3z^{64}+4z^{65}+z^{66}
$$
(the permutation $x1,x2,x3,t,x4,x5,x6,x7,x8,x9$      still gives a finite Gr\"obner basis),
and the Ufnarovskij method now still gives that $1/A(z^{31},z)=1+3z^{64}+4z^{65}+z^{66}-5z^{32}-4z^{33}$.
On the other hand the Hilbert series of $HKR(x,y)=1+5xy+4xy^2+3x^2y^2+4x^2y^3+x^2y^4$
so that $HKR(-z^{31},z)=1-5z^{32}-4z^{33}+3z^{64}+4z^{65}+z^{66}$
and the formula (10) now also gives
$$
1/\Phi_R(z^{31},z)=(1+1/z^{31})/A(z^{31},z)-HKR(-z^{31},z)/z^{31}=1+3z^{64}+4z^{65}+z^{66}-5z^{32}-4z^{33}
$$
But the taylor series of $1/(1+3z^{64}+4z^{65}+z^{66}-5z^{32}-4z^{33})$ is
$$
1+z^{32}(5+4z)+z^{64}(22+36z+15z^2)+z^{96}(95+236z+198z^{2}+56z^{3})+
$$
$$
+z^{128}(409+1364z+1723z^{2}+976z^{3}+109z^{4})+z^{160}(1760+7368z+12438z^{2}+\ldots)+\ldots
$$
which ``explains'' the columns above at the beginning of this section
and gives an indication that 
$$
1/\Phi_R(x,y)=1-(5+4y+y^2)xy+(3+4y+y^2)x^2y^2
$$
In a similar way one treats all the quadratic embedding dimension $4$ cases and the results
are given in the tables A-D that follow.
\vfill \break
\newdimen\tempdim                
\newdimen\othick \othick=.4pt    
\newdimen\ithick \ithick=.4pt    
\newdimen\spacing \spacing=9pt   
\newdimen\abovehr \abovehr=6pt   
\newdimen\belowhr \belowhr=8pt   
\newdimen\nexttovr \nexttovr=8pt 

\def\r{\hfil&\omit\vrsp\vrule width\othick\cr&}   
\def\rr{\hfil\down{\abovehr}&\omit\vrsp\vrule width\othick\cr
 \noalign{\hrule height\ithick}\up{\belowhr}&}     
\def\up#1{\tempdim=#1\advance\tempdim by1ex
  \vrule height\tempdim width0pt depth0pt} 
\def\down#1{\vrule height0pt depth#1 width0pt} 
\def\large#1#2{\setbox0=\vtop{\hsize#1 \lineskiplimit=0pt \lineskip=1pt
  \baselineskip\spacing \advance\baselineskip by 3pt \noindent
  #2}\tempdim=\dp0\advance\tempdim by\abovehr\box0\down{\tempdim}}
\def\vrsp{\hskip\nexttovr\relax}
\def\toprule#1{\def\startrule{\hrule height#1\relax}} 
\toprule{\othick}     
\def\nstrut{\vrule height\spacing depth3.5pt width0pt}
\def\exclaim{\char`\!}   
\def\preamble#1{\def\startup{#1}}  
\preamble{&##}  
{\catcode`\!=\active
\gdef!{\hfil\vrule width0pt\vrsp\vrule width\ithick\relax\vrsp&}}

\def\table #1{\vbox\bgroup \setbox0=\hbox{#1}
 \vbox\bgroup\offinterlineskip \catcode`\!=\active
\halign\bgroup##\vrule width\othick\vrsp&\span\startup\nstrut\cr
\noalign{\medskip}
\noalign{\startrule}\up{\belowhr}&}

\def\caption #1{\down{\abovehr}&\omit\vrsp\vrule width\othick\cr
\noalign{\hrule height\othick}\egroup\egroup \setbox1=\lastbox
\tempdim=\wd1 \hbox to\tempdim{\hfil \box0 \hfil} \box1 \smallskip
\hbox to\tempdim{\advance\tempdim by-20pt\hfil\vbox{\hsize\tempdim
\noindent #1}\hfil}\egroup}
$$\table{Table A. Double Tor of HKR }
R !  ideal gen:s (ex.)  !1/Double Tor of $HKR$ \rr
\hfill{ 1}!$0$!$1$\rr
\hfill{ 2}!$x^2$!$1-xy$\rr
\hfill{ 3}!$x^2,y^2$!$1-2xy+x^2y^2$\rr
\hfill{ 4}!$x^2,xy$!$1-(2+y)xy$\rr
\hfill{ 5}!$x^2,y^2,z^2$!$1-3xy+3x^2y^2-x^3y^3$\rr
\hfill{ 6}!$x^2,y^2+xz,yz$!$1-(3+y+2y^2)xy+(3-y)x^2y^2-x^3y^3$\rr
\hfill{\bf 7}!$x^2+y^2,z^2+u^2,xz+yu$!$1-(3+2y+4y^2+y^3)xy+(3-y)x^2y^2-x^3y^3$\rr
\hfill{ 8}!$x^2,y^2,xz$!$1-(3+y)xy+(2+y)x^2y^2$\rr
\hfill{ 9}!$x^2,xy,y^2$!$1-(3+2y)xy$\rr
\hfill{ 10}!$x^2,xy,xz$!$1-(3+3y+y^2)xy$\rr
\hfill{ 11}!$x^2,y^2,z^2,u^2$!$1-4xy+6x^2y^2-4x^3y^3+x^4y^4$\rr
\hfill{\bf 12}!${\scriptstyle x^2+xy,y^2+yu,x^2+xu,u^2+zu}$!$1-(4+y+3y^2+3y^3)xy+(6-y^2)x^2y^2-4x^3y^3+x^4y^4$\rr
\hfill{\bf 13}!${\scriptstyle x^2+z^2+u^2,y^2,xz,yu+zu}$!$1-(4+y+2y^2)xy+(6+2y^2)x^2y^2-(4-y)x^3y^3+x^4y^4$\rr
\hfill{\bf 14}!$xz,y^2,z^2+u^2,yu+zu$!$1-(4+2y+5y^2)xy+(6-3y+2y^2)x^2y^2-(4-y)x^3y^3+x^4y^4$\rr
\hfill{\bf 15}!$xz,y^2,yz+u^2,yu+zu$!$1-(4+3y+8y^2+2y^3)xy+(6-4y)x^2y^2-(4-y)x^3y^3+x^4y^4$\rr
\hfill{\bf 16}!${\scriptstyle xy+z^2+yu,y^2,yu+zu,xz}$!$1-(4+4y+10y^2+3y^3)xy+(6-4y)x^2y^2-(4-y)x^3y^3+x^4y^4$\rr
\hfill{\bf 17}!$xz,yz+xu,y^2,yu+zu$!${\scriptstyle({1-(4+y^2)xy+(6-4y-6y^2-5y^3-y^4)x^2y^2-(4-y)x^3y^3+x^4y^4})/(1+xy^2)}$\rr
\hfill{\bf 18}!$x^2,y^2,z^2,yu$!$1-(4+y)xy+(5+2y)x^2y^2-(2+y)x^3y^3$\rr
\hfill{\bf 19}!$xz,y^2,yu+zu,u^2$!$1-(4+2y+3y^2+y^3)xy+(5-y^2)x^2y^2-(2+y)x^3y^3$\rr
\hfill{\bf 20}!$xz,y^2,yu+z^2,yu+zu$!$1-(4+3y+5y^2+2y^3)xy+(5+2y)x^2y^2-(2+y)x^3y^3$\rr
\hfill{\bf 21}!$xz,y^2,z^2,yu+zu$!$1-(4+2y)xy+(4+4y+y^2)x^2y^2$\rr
\hfill{\bf 22}!$x^2+xy,xu,xz+yu,y^2$!${\scriptstyle 1-(4+y)xy+(4+4y-2y^3-y^4)x^2y^2-(1+4y+y^2)x^3y^4+(-1-2y+y^2)x^4y^4}$\rr
\hfill{ 23}!$xz,xu,y^2,z^2$!$1-(4+2y)xy+(3+2y)x^2y^2$\rr
\hfill{\bf 24}!$xz,y^2,yz+z^2,yu+zu$!$1-(4+y)xy+(5+2y)x^2y^2-(2+y)x^3y^3$\rr
\hfill{\bf 25}!$x^2,xy,xz,u^2$!$1-(4+3y+y^2)xy+(3+3y+y^2)x^2y^2$\rr
\hfill{ 26}   !$xz,y^2,yu,zu$!$1-(4+3y)xy+(1+y)x^2y^2$\rr
\hfill{ 27}   !$xy,xz,y^2,yz$!$1-(4+4y+y^2)xy$\rr
\hfill{\bf 28}!$x^2,xy,xz,xu$!$1-(4+y)xy+(5+2y)x^2y^2-(2+y)x^3y^3$\rr
\hfill{\bf 29}!${\scriptstyle x^2+xy,y^2+xu,z^2+xu,zu+u^2}$!${\scriptstyle 1-(5+5y+16y^2+5y^3)xy+(10-10y)x^2y^2-(10-5y)x^3y^3+(5-y)x^4y
^4-x^5y^5}$\rr
\hfill{\bf 30}!${\scriptscriptstyle xy+u^2,xz,x^2+z^2+u^2,y^2,yu+zu}$!${\scriptscriptstyle{1-(5+y+5y^2)xy+(10-10y-10y^2-10y^3-y^4)x^2y^
2-(10-5y-y^3)x^3y^3+(5-y)x^4y^4-x^5y^5}\over{1+xy^2}}$\rr
\hfill{\bf 31}!${\scriptstyle x^2-y^2,y^2-z^2,z^2-u^2,xz+yu,
                }$!${\scriptscriptstyle{1-(5+2y+8y^2+3y^3)xy+(10-10y-12y^2-13y^3-4y^4)x^2y^2-(10-5y)x^3y^3+(5-y)x^4y^4-x^5y^5}\over{1+x
y^2}}$\r
! ${\scriptstyle -x^2+xy-yz+xu}$!\rr
\hfill{\bf 32}!$x^2+z^2,xz,y^2,yu+zu,u^2$!${\scriptstyle 1-(5+3y+7y^2+2y^3)xy+(9-3y-2y^2)x^2y^2-(7+y-y^2)x^3y^3+(2+y)x^4y^4}$
\caption{\sl }
$$
$$\table{Table B. Double Tor of HKR }
R !  ideal gen:s (ex.)  !1/Double Tor of $HKR$ \rr
\hfill{\bf 33}!$x^2+xy,y^2+yz,$!${\scriptstyle 1-(5+4y+10y^2+4y^3)xy+(9-4y-3y^2)x^2y^2-(7+y-y^2)x^3y^3+(2+y)x^4y^4}$\r
!$y^2+xu,z^2+xu,zu+u^2$!\rr
\hfill{\bf 34}!$x^2+xy+yu+u^2,y^2,xz$!${\scriptstyle (1-5xy+(9-5y-18y^2-15y^3-4y^4)x^2y^2-(5+5y+9y^3+10y^4+3y^5)x^3y^3+}$\r
! $ x^2+z^2+u^2,yu+zu$!${\scriptstyle+(-5+9y)x^4y^4+(9-5y)x^5y^5-(5-y)x^6y^6+x^7y^7)/(1+xy^2)^2}$\rr
\hfill{\bf 35}!$x^2+z^2+u^2,y^2,xz$!${\scriptstyle (1-(5+y+2y^2+y^3)xy+(9-4y-12y^2-10y^3-3y^4)x^2y^2-(7+y-y^2)x^3y^3+}$\r
!$xy+yz+yu,yu+zu$!${\scriptstyle+(2+y)x^4y^4)/(1+xy^2)}$\rr
\hfill{\bf 36}!$x^2+y^2,z^2,u^2$!${\scriptstyle 1-(5+3y+4y^2+2y^3)xy+(8+2y-3y^2-y^3)x^2y^2-(4+4y+y^2)x^3y^3}$\r
!$yz-yu,xz+zu$!\rr
\hfill{\bf 37}!$x^2,y^2,xy-zu$!${\scriptstyle (1-(5+2y+4y^2+2y^3)xy+(8-3y-7y^2-7y^3-3y^4)x^2y^2-(4-2y+2y^2+3y^3)x^3y^3 +}$\r
!$yz-xu,(x-y)(z-u)$!${\scriptstyle +(-1-4y-y^2+y^3)x^4y^4+(1+2y+y^2)x^5y^5)/(1+xy^2)}$\rr
\hfill{\bf 38}!$x^2,y^2,z^2,zu,u^2$!$1-(5+2y)xy+(7+4y)x^2y^2-(3+2y)x^3y^3$\rr
\hfill{\bf 39}!$x^2+yz+u^2,xz+z^2+yu,$!$1-(5+3y+3y²+y³)xy+(7+2y-y²)x^2y^2-(3+2y)x^3y^3$\r
!$xy,xu,zu$!\rr
\hfill{\bf 40}!$x^2-xu,xu-y^2,y^2-z^2$!$1-(5+4y+6y^2+3y^3)xy+(7+y-2y^2)x^2y^2-(3+2y)x^3y^3$\r
!$z^2-u^2,xz-yu$!\rr
\hfill{\bf 41}!$xy,y^2,z^2,zu,u^2$!$1-(5+3y)xy+(6+7y+2y^2)x^2y^2$\rr
\hfill{\bf 42}!$x^2+xy,zu,y^2,$!$(1-(5+2y)xy+(6+y-5y^2-5y^3-2y^4)x^2y^2+$\r
!$xu,xz+yu$!$+(-1+y-y^3)x^3y^3+(-1-2y-y^2)x^4y^4)/(1+xy^2)$\rr
\hfill{\bf 43}!$x^2,y^2,yz,zu,u^2$!$1-(5+3y)xy+(5+4y)x^2y^2-(1+y)x^3y^3$\rr
\hfill{\bf 44}!$xz,yz,y^2$!$1-(5+4y+3y^2+2y^3)xy+(5+3y-y^2)x^2y^2-(1+y)x^3y^3$\r
!$yu+zu,z^2+u^2$!\rr
\hfill{\bf 45}!$xy+yz,xy+z^2+yu,$!$1-(5+4y)xy+(4+6y+2y^2)x^2y^2$\r
!$yu+zu,y^2,xz$!\rr
\hfill{\bf 46}!$x^2,xy,yz,zu,u^2$!$1-(5+4y)xy+(3+4y+y^2)x^2y^2$\rr
\hfill{\bf 46}!$xz+u^2,xy,xu,x^2,$!$1-(5+4y+y^2)xy+(4+4y+y^2)x^2y^2$\r
\hfill{\bf va}!$zu+y^2+z^2$!\rr
\hfill{\bf 47}!$x^2+xy,y^2,xu,$!$1-(5+5y+2y^2+y^3)xy+(3+4y+y^2)x^2y^2$\r
!$xz+yu,-x^2+xz-yz$!\rr
\hfill{\bf 48}!$xy,z^2+yu,yu+zu,$!$1-(5+5y+3y^2+2y^3)xy+(3+2y-y^2)x^2y^2-(1+y)x^3y^3$\r
!$y^2,xz$!\rr
\hfill{\bf 49}!$xz,y^2,z^2,yu,zu$!$1-(5+5y+y^2)xy+(2+3y+y^2)x^2y^2$\rr
\hfill{ 50}!$x^2,xy,xz,y^2,z^2$!$1-(5+5y+y^2)xy+(1+y)x^2y^2$\rr
\hfill{\bf 51}!$xy,xz,yz+xu,z^2,zu$!$1-(5+6y+3y^2+y^3)xy+(1+y)x^2y^2$\rr
\hfill{ 52}!$x^2,xy,xz,y^2,yz$!$1-(5+6y+2y^2)xy$\rr
\hfill{\bf 53}!$y^2-u^2,xz,yz,z^2,zu$!$1-(5+7y+4y^2+y^3)xy$\rr
\hfill{\bf 54}!$x^2,xz,y^2,z^2,yu+zu,u^2$!$1-(6+4y)xy+(9+12y+4y^2)x^2y^2$\rr
\hfill{\bf 55}!$x^2+xy,xz+yu,xu,$!${\scriptscriptstyle(1-(6+2y)xy+(9-y-9y^2-9y^3-4y^4)x^2y^2-(1-6y+7y^2+18y^3+12y^4+4y^5)x^3y^3+}$\r
!$y^2,z^2,zu+u^2$!${\scriptscriptstyle+(-6-13y-3y^2+3y^3)x^4y^4+(3+6y+3y^2)x^5y^5)/(1+xy^2)^2}$
\caption{\sl }
$$
$$\table{Table C. Double Tor of HKR }
R !  ideal gen:s (ex.)  !1/Double Tor of $HKR$ \rr
\hfill{\bf 56}!$x^2+xz+u^2,xy,xu,$!${\scriptstyle(1-(6+3y+y^2+y^3)xy+(9+y-11y^2-10y^3-4y^4)x^2y^2+}$\r
!$x^2-y^2,z^2,zu$!${\scriptstyle+(-4-y+y^2-y^3)x^3y^3-(1+y)x^4y^5))/(1+xy^2)}$\rr
\hfill{\bf 57}!${\scriptstyle x^2+yz+u^2,xy,xu,zu,y^2+z^2,xz+yu}$!${\scriptstyle 1-(6+5y+4y^2+3y^3)xy+(8+5y-3y^2-y^3)x^2y^2-(3+4y+y^2\
)z^3y^3}$\rr
\hfill{\bf 57}!$ x^2+y^2+z^2,xy,$!${\scriptstyle 1-(6+5y+5y^2+2y^3)xy+(9+4y-2y^2-y^3)x^2y^2-(4+4y+y^2)z^3y^3}$\r
\hfill{\bf v2}!$xu,yz,zu,xz+u^2$!\rr
\hfill{\bf 58}!${\scriptstyle xy,x^2+zu,y^2,z^2,yu+xz,xu}$!$ 1-(6+5y)xy+(6+9y+3y^2)x^2y^2$\rr
\hfill{\bf 59}!$x^2-y^2,xy,xu,$!${\scriptstyle(1-(6+4y)xy+(6+3y-3y^2-3y^3-2y^4)x^2y^2+}$\r
!$z^2,zu,xz+yu$!${\scriptstyle+(4+2y-3y^2-y^3)x^3y^4-(1+4y+4y^2+y^3)x^4y^4)/(1+xy^2)}$\rr
\hfill{\bf 59}!$x^2-y^2,xy,yz,$!${\scriptstyle(1-(6+4y+y^2)xy+(7+3y-3y^2-3y^3-y^4)x^2y^2+}$\r
\hfill{\bf va}!$zu,xz+u^2,xu$!${\scriptstyle+(4+2y-2y^2-y^3)x^3y^4-(2+5y+4y^2+y^3)x^4y^4)/(1+xy^2)}$\rr
\hfill{\bf 60}!$x^2+yz+u^2,xy,zu,$!${\scriptstyle(1-(6+4y)xy+(6+2y-5y^2-4y^3-2y^4)x^2y^2+}$\r
\hfill{\bf va}!$z^2,xz+yu,xu$!${\scriptstyle+(-1+3y+4y^2-y^3-y^4)x^3y^3-(1+2y+y^2)x^4y^5)/(1+xy^2)}$\rr
\hfill{\bf 61}!$x^2-y^2,xy,z^2,xu,zu,u^2$!$1-(6+5y+y^2)xy+(6+6y+y^2)x^2y^2-(y+1)x^3y^3$\rr
\hfill{\bf 62}!$x^2-y^2,xy,xu,yz+yu,z^2,zu$!$1-(6+6y+3y^2+3y^3)xy+(5+6y-y^2)x^2y^2-x^3y^4$\rr
\hfill{\bf 62}!$ x^2+yz+u^2,yu,$!$1-(6+6y+4y^2+2y^3)xy+(6+5y)x^2y^2-(1+y)x^3y^3$\r
\hfill{\bf va}!$ zu,xy,z^2,xu$!\rr
\hfill{\bf 63}!$x^2,xy,xu,y^2,z^2,zu$!$1-(6+6y+y^2)xy+(4+6y+2y^2)x^2y^2$\rr
\hfill{\bf 63}!$y^2,xz+yu,zu,xy,z^2,xu$!$1-(6+6y)xy+(3+6y+2y^2)x^2y^2$\r
\hfill{\bf v4}!!\rr
\hfill{\bf 63}!$x^2,xy,xu,yu,z^2,$!$1-(6+6y+2y^2)xy+(5+6y+2y^2)x^2y^2$\r
\hfill{\bf v8}!$xz+u^2,y^2+z^2+zu$!\rr
\hfill{\bf 63}!$x^2,xy,xz+u^2,xu,$!$(1-(6+5y+y^2)xy+(4-4y^2-2y^3-y^4)x^2y^2+$\r
\hfill{\bf ne}!$y^2+z^2,zu$!$+(4+5y+y^2)x^3y^4)/(1+xy^2)$\rr
\hfill{\bf 64}!$x^2-y^2,xy,z^2,xu,yu,zu$!$1-(6+6y+y^2)xy+(3+3y)x^2y^2-(1+y)x^3y^3$\rr
\hfill{\bf 65}!$x^2,xy,xz,y^2,yu+z^2,yu+zu$!$1-(6+7y+4y^2+2y^3)xy+$\r
\hfill{\bf }!$x^2,xy,xz,y^2,yu+z^2,yu+zu$!$+(3+2y-y^2)x^2y^2-(1+y)x^3y^3$\rr
\hfill{\bf 66}!$xz,y^2,yu,z^2,zu,u^2$!$1-(6+7y+2y^2)xy+(2+3y+y^2)x^2y^2$\rr
\hfill{\bf 66}!$xy,xz+u^2,xu,yu,zu,z^2$!$1-(6+7y+y^2)xy+(1+3y+y^2)x^2y^2$\r
\hfill{\bf v5}! !\rr
\hfill{\bf 67}!$xy,xz,y^2,yu,z^2,zu$!$1-(6+8y+3y^2)xy+(1+2y+y^2)x^2y^2$\rr
\hfill{ 68}!$x^2,xy,xz,y^2,yz,z^2$!$1-(6+8y+3y^2)xy$\rr
\hfill{\bf 68v}!$x^2,xy,xz,xu,u^2,y^2+z^2+zu$!$1-(6+8y+4y^2+y^3)xy+(1+y)x^2y^2$\rr
\hfill{\bf 69}!$x^2,xz,xu,xy-zu,yz,z^2$!$1-(6+9y+5y^2+y^3)xy$\rr
\hfill{\bf 70}!$x^2,xy,xz,xu,y^2,yz$!$1-(6+9y+5y^2+y^3)xy$\rr
\hfill{\bf 71}!$x^2,y^2,z^2,u^2,xy,zu,yz+xu$!$1-(7+8y+2y^2)xy+(5+8y+3y^2)x^2y^2$
\caption{\sl }
$$
$$\table{Table D. Double Tor of HKR }
R !  ideal gen:s (ex.)  !1/Double Tor of $HKR$ \rr
\hfill{\bf 71}!$x^2,y^2+z^2,xy,yz,$!$(1-(7+7y+2y^2)xy+(5+y-3y^3-y^4)x^2y^2+$\r
\hfill{\bf v16}!$zu,xz+u^2,xu$!$+(5+7y+2y^2)x^3y^4)/(1+xy^2)$\rr
\hfill{\bf 71}!$x^2+u^2,xy,xu,y^2,yz,z^2,zu$!$1-(7+8y+3y^2)xy+(6+8y+3y^3)x^2y^2$\r
\hfill{\bf v4}! !\rr
\hfill{\bf 71}!$x^2,y^2,z^2,xz+u^2,xu,yz,zu$!$1-(7+8y+y^2)xy+(4+8y+3y^3)x^2y^2$\r
\hfill{\bf v7}! !\rr
\hfill{\bf 71}!$x^2+xy,x^2+yz,xy+y^2,z^2,$!$1-(7+8y)xy+(3+8y+3y^3)x^2y^2$\r
\hfill{\bf v5}!$z^2,xu,zu,xz+u^2$!\rr
\hfill{\bf 72}!$x^2-y^2,z^2,xy,yz,zu,xu$!$1-(7+9y+2y^2)xy+(2+5y+2y^3)x^2y^2$\rr
\hfill{\bf 72}!$xu+u^2,x^2+xy,y^2+xu,y^2+yz,$!$1-(7+9y+3y^2)xy+(3+5y+2y^3)x^2y^2$\r
\hfill{\bf v1}!$y^2+yz,yu+zu,z^2+xu,zu+u^2$!\rr
\hfill{\bf 72}!$yz,x^2+xy,xz+yu,xu,z^2,$!$1-(7+9y+1y^2)xy+(1+5y+2y^3)x^2y^2$\r
\hfill{\bf v2e}!$zu,x^2+u^2$!\rr
\hfill{\bf 73}!$x^2,y^2,z^2,u^2,xu,yu,zu$!$1-(7+9y+4y^2+y^3)xy+$\r
\hfill{\bf }!!$+(3+3y)x^2y^2-(1+y)x^3y^3$\rr
\hfill{\bf 74}!$x^2,xy+z^2,yz,xu,yu,zu,u^2$!$1-(7+10y+7y^2+3y^3)xy+$\r
\hfill{\bf }!!$+(3+2y-y^2)x^2y^2-(1+y)x^3y^3$\rr
\hfill{\bf 75}!$x^2,xy,xz,y^2,yz,yu,u^2$!$1-(7+10y+5y^2+y^3)xy+(2+3y+y^2)x^2y^2$\rr
\hfill{\bf 75}!$y^2,xz,yz+xu,z^2,yu,zu,u^2$!$1-(7+10y+4y^2)xy+(1+2y+y^2)x^2y^2$\r
\hfill{\bf v1}!!\rr
\hfill{\bf 75}!$xy,xz+u^2,xu,yz+u^2,yu,zu,z^2$!$1-(7+10y+3y^2)xy+(2y+y^2)x^2y^2$\r
\hfill{\bf v2}!!\rr
\hfill{\bf 76}!$x^2,xy,xz,xu,z^2,zu,yu$!$1-(7+11y+6y^2+y^3)xy$\rr
\hfill{\bf 77}!$x^2,xy,xz,xu,y^2,yz,yu$!$1-(7+12y+8y^2+2y^3)xy$\rr
\hfill{\bf 78}!$x^2,xy,y^2,z^2,zu,u^2,xz+yu,yz-xu$!$1-(8+12y+5y^2)xy+(2+4y+y^2)x^2y^2$\rr
\hfill{\bf 78}!$x^2,y^2,z^2,u^2,xy,xz,xu,yu$!$1-(8+12y+6y^2+y^3)xy+(3+5y+2y^2)x^2y^2$\r
\hfill{\bf v1}!!\rr
\hfill{\bf 78}!$xu,yu+xz,yz,x^2,y^2,z^2+xy,u^2,zu$!$1-(8+12y+3y^2)xy+(4y+2y^2)x^2y^2$\r
\hfill{\bf v2e}!!\rr
\hfill{\bf 78}!$xu,yu+xz,yz,x^2,y^2,z^2+xz,u^2,zu$!$1-(8+12y+4y^2)xy+(1+4y+2y^2)x^2y^2$\r
\hfill{\bf v3v}!!\rr
\hfill{\bf 79}!$x^2,xy,xz,xu,y^2,yu,z^2,zu$!$1-(8+13y+7y^2+y^3)xy+(1+2y+y^2)x^2y^2$\rr
\hfill{\bf 80}!$x^2,xy,xz,yz,y^2,yu,z^2,zu$!$1-(8+14y+9y^2+2y^3)xy$\rr
\hfill{\bf 81}!$x^2,y^2,z^2,u^2,xy,xz,yz-xu,yu,zu$!$1-(9+16y+9y^2)xy+x^2y^4$\rr
\hfill{\bf 81}!$x^2,y^2,z^2,u^2,xy,xz,xu,yu,zu$!$1-(9+16y+10y^2+2y^3)xy+(1+2y+y^2)x^2y^2$\r
\hfill{\bf va}!!\rr
\hfill{\bf 82}!$x^2,xy,xz,xu,y^2,zu,u^2,yz,yu$!$1-(9+17y+12y^2+3y^3)xy$\rr
\hfill{\bf 83}!$x^2,y^2,z^2,u^2,xy,xz,xu,yz,yu,zu$!$1-(10+20y+15y^2+4y^3)xy$
\caption{\sl Case 81 is Gorenstein}
$$

\mysec{6. Higher embedding dimensions. Three examples.}

Many new phenomena come up in connection with higher embedding dimensions of local rings.
Here we will just mention three very different cases of 
$R=k[x,y,z,u,v]/J$ where the $J$ is an ideal generated by quadratic forms in $x,y,z,u,v$

{\bf Case I}  $R_I=k[x,y,z,u,v]/(y^2,u^2,yz+xu,zu+yv,z^2-yu-xv)$

{\bf Case II} $R_{II}=k[x,y,z,u,v]/(xu-xv,xv-zu,yv,u^2,v^2)$

{\bf Case III} $R_{III}=k[x,y,z,u,v]/(xz+yu,yv,zu+uv,z^2,v^2)$

In all these three cases the resolution of J over $k[x,y,z,u,v]$ has the same diagram of Betti numbers:
$$
\left ( \matrix{
 1 & .&.&.&.&.\cr
 .&5&3
&.&.&. \cr
 .&.&8 &12&6&1 \cr} \right )
$$
But the homological properties of $HKR$ and $R$ are {\it very} different:

{\bf Case I} This is the case that L\"ofwall and I studied in [LR2], which also solved an old problem
(the ``Opera case''). In [LR2] we proved that the Hilbert series of $R^!$ is transcendental and different for all values
of the characteristic of the field $k$.
Similar results are true for $HKR^!$ and the Avramov spectral sequence does not degenerate.
We will just give a few indications below.

{\bf Case II} 
In [R1] I found, inspired by Lemaire [LE], a case of a local ring $(R,m)$ where the Yoneda Ext-algebra ${\rm Ext}^*_R(k,k)$
was not finitely generated. Furthermore that $R^!$ had global dimension 3.
This is a variant of that example and similar assertions are true for $HKR$.
Furthermore the Avramov spectral sequence degenerates.

{\bf Case III} This is a variant of the first known example due to Anick [A] (cf also [LR1]) where $P_R(z)$ is transcendental.
Here $HKR$ seems to have the same properties and the Avramov spectral sequence degenerates at least up to degrees 20.

More details:
{\bf Case I}. 

We apply as before {\tt DGAlgebras} and obtain that $HKR$ has the cube of the augmmentation ideal equal to 0.
Furthermore $HKR$ has 9 generators: five $X_1,X_2,X_3,X_4,X_5$ of degree (1,2), three generators
$X_6,X_7,X_8$ of degree (2,3) and one $X_9$ of degree (5,7). Furthermore these generators only
satisfy 17 quadratic relations of which the following 9 are non-monomial relations:
$$
X_3X_4-X_2X_5,X_1X_4-X_3X_5,X_2X_3+X_1X_5,X_4X_7-X_5X_8,X_5X_6+X_3X_8,
$$
$$
X_4X_6-X_5X_7+X_2X_8,X_3X_6+X_1X_7,X_2X_6+X_3X_7+X_1X_8,X_7X_7+2X_6X_8
$$
and since the following 11 quadratic monomials do {\it not} occur among the 17 quadratic relations
$$
X_1X_2,X_1X_3,X_2X_4,X_4X_5,X_1X_6,X_2X_7,X_4X_8,X_6X_6,X_6X_7,X_7X_8,X_8X_8
$$
it follows that $HKR^!$ is the quotient of the free algebra on the 9 generators $Y_i$ (which are
dual to the $X_i$) by the two-sided ideal generated by the 21 quadratic relations:
$$
[Y_1,Y_2],[Y_1,Y_3],[Y_2,Y_4],[Y_4,Y_5],[Y_1,Y_6],[Y_2,Y_7],[Y_4,Y_8],Y_6Y_6,[Y_6,Y_7],[Y_7,Y_8],Y_8Y_8,
$$
$$
[Y_1,Y_4]+[Y_3,Y_5],[Y_2,Y_3]-[Y_1,Y_5],[Y_3,Y_4]+[Y_2,Y_5],[Y_2,Y_6]-[Y_3,Y_7],
$$
$$
[Y_3,Y_7]-[Y_1,Y_8],[Y_3,Y_6]-[Y_1,Y_7],[Y_4,Y_6]+[Y_5,Y_7],[Y_5,Y_7]+[Y_2,Y_8],
$$
$$
[Y_5,Y_6]-[Y_3,Y_8],[Y_4,Y_7]+[Y_5,Y_8],[Y_7,Y_7]-[Y_6,Y_8]
$$
Here the graded commutators are as before and adding as before a new variable T giving a new algebra
$A$ generated in degree 1 such that
$$
{1\over HKR^!(z)}= {1\over A(z)}+10z-5z^2-3z^3-z^6
$$
We find using {\tt BERGMAN} that
$$
{1\over A(z)}=1-10z+7z^4+10z^5+5z^6-10z^8-27z^9-25z^{10}+3z^{11}+19z^{12}+27z^{13}+56z^{14}+57z^{15}+
$$
$$
+z^{16}-42z^{17}-58z^{18}-101z^{19}-124z^{20}-51z^{21}+48z^{22}+97z^{23}+154z^{24}+204z^{25}+\ldots
$$
This gives ${1\over HKR^!(z)}$ and since the two-variable Hilbert series of $HKR$ is
$1+5xy+3xy^2+xy^5+7x^2y^2+10x^2y^3+5x^2y^4$ and therefore the one variable series
with a minus sign in the first variable is
 $h=1-5z^2-3z^3+7z^4+10z^5+4z^6$,
the formula
$$
{1\over P_{HKR}(z)}=(1+1/z)/HKR^!(z)-h/z
$$
gives
$$
P_{HKR}(z)=1+5z^2+3z^3+18z^4+20z^5+60z^6+93z^7+221z^8+415z^9+929z^{10}+1936z^{11}+\ldots
$$
whereas
$$
{P_R(z)\over(1+z)^5}=1+5z^2+3z^3+18z^4+20z^5+59z^6+92z^7+216z^8+407z^9+907z^{10}+1897z^{11}+\ldots
$$
and they start differing in degree 6 so that the Avramov spectral sequence does not degenerate.
The preceding calculations are in characteristic 0. But e.g. in characteristic 7 we have 
$$
{P_R(z)\over(1+z)^5}=1+5z^2+3z^3+18z^4+20z^5+60z^6+93z^7+221z^8+415z^9+928z^{10}+1935z^{11}+\ldots
$$
In characteristic 7 we also have
$$
P_{HKR}(z)=1+5z^2+3z^3+18z^4+20z^5+60z^6+93z^7+221z^8+415z^9+929z^{10}+1936z^{11}+\ldots
$$
and these two series differ in degree 11, so that the Avramov spectral sequence does not degenerate.

Similar results are true in characteristic 2,3,5,11,13, etc. But characteristic 2 is slightly different. 

Important observation: The relations between the first 5 variables $X_i$ in $HKR$ are the 
same in all
characteristic:
$$
X_3X_4-X_2X_5,X_1X_4-X_3X_5,X_2X_3+X_1X_5
$$
But the quotient of the exterior algebra in five variables with these 3 relations has a homology
with transcendental Hilbert series which is different for all characteristic
(it is the ``skew'' counterpart of the present {\bf Case I}  for ordinary polyomial rings).
It occurs among the relations that Eisenbud and Koh have classified in [Ei-K].

More details:
{\bf Case II}.
Now the ring is $R=k[x,y,z,u,v]/(xu-xv,xv-zu,yv,u^2,v^2)$. We have that the Hilbert series $R^!(z)={(1+z)^3(1-z)\over(1-z-z^2)^3}$ and $R(z)={1+2z-2z^2-3z^3+4z^4-z^5\over(1-z)^3}$
and 
$$
1/P_R(z)=(1+1/z)/R^!(z)-R(-z)/z = {(1-2z-z^2)(1-z-z^2+z^3-z^4)\over (1+z)^3(1-z)}
$$
We now apply {\tt DGAlgebras} and we obtain that $HKR$ has five generators $X_1,X_2,X_3,X_4,X_5$
of degree (1,2),
three generators $X_6,X_7,X_8$ of degree (2,3), one generator $X_9$ of degree (2,4), three generators
$X_{10},X_{11},X_{12}$  of degree (3,5),
three generators $X_{13},X_{14},X_{15}$ of degree (4,6) and one generator $X_{16}$ of degree (5,7).
Furthermore these generators satisfy 107  quadratic relations of which only the following four are
non-monomial relations:
$$
X_1X_3+X_3X_4,X_1X_2+X_3X_4,X_1X_8-X_4X_8,X_2X_6-X_3X_6
$$
and since the following 16 quadratic monomials do {\it not} occur among the 107 quadratic relations:
$$
X_1X_4,X_1X_5,X_2X_3,X_2X_4,X_2X_5,X_3X_5,X_1X_7,X_2X_7,X_2X_8,X_3X_7,
$$
$$
X_4X_6,X_5X_6,X_5X_8,X_6X_7,X_6X_8,X_7X_8
$$
it follows that $HKR^!$ is the quotient of the free algebra on the
16 generators $Y_i$ (which are dual to the $X_i$) by the two-sided ideal generated by
the 19 quadratic relations
$$
[Y_1,Y_4],[Y_1,Y_5],[Y_2,Y_3],[Y_2,Y_4],[Y_2,Y_5],[Y_3,Y_5],[Y_1,Y_7],[Y_2,Y_7],[Y_2,Y_8],[Y_3,Y_7],
$$
$$
[Y_4,Y_6],[Y_5,Y_6],[Y_5,Y_8],[Y_6,Y_7],[Y_6,Y_8],[Y_7,Y_8],[Y_3,Y_4]-[Y_1,Y_3]-[Y_1,Y_2], \leqno(xxx)
$$
$$
[Y_1,Y_8]+[Y_4,Y_8],[Y_2,Y_6]+[Y_3,Y_6]
$$
Here $[,]$ means the graded commutator so that $[Y_i,Y_j]=Y_iY_j-Y_jY_i$ if $1\leq i \leq 5$ or $i=10,11,12,16$ {\it and} $1\leq j \leq 5$ or $j=10,11,12,16$. Furthermore we have $[Y_s,Y_t]=Y_sY_t+Y_tY_s$
otherwise. 
To calculate $HKR^!(z)$ we use again the Ufnarovski trick and replace the $Y_i$ with $Z_iT^{deg Y_i-1}$
everywhere in order to get a new algebra $A$ with generators of degree 1. Calculations in {\tt BERGMAN}
gives after a few minutes
$$
{1\over A(z)}=1-17z+7z^4+9z^5+2z^6-3z^7-4z^8-4z^9-4z^{10}-4z^{11}-4z^{12}-4z^{13}-4z^{14}+O(z^{15})
$$
leading as before to
$$
{1\over HKR^!(z)}={1\over A(z)}+17z-5z^2-4z^3-3z^4-3z^5-z^6={(1+z)^2(1-3z+4z^3-2z^5-z^6)\over 1-z}
$$

 Furthermore
the cube of the augmentation ideal of $HKR$ is equal to zero, and the square
of that ideal is generated by 
$$
X_3X_5,X_2X_5,X_1X_5,X_3X_4,X_2X_4,X_1X_4,X_2X_3,X_5X_8,X_4X_8,X_2X_8,X_3X_7,X_2X_7,X_1X_4,
$$
$$
X_5X_6,X_4X_6,X_3X_6,X_7X_8,X_6X_8,X_6X_7
$$
leading to the bigraded Hilbert series
 $$
HKR(x,y)=1+5xy+4xy^2+3xy^3+3xy^4+xy^5+7x^2y^2+9x^2y^3+3x^2y^4
$$
and thus $H(-z,z)=1-5z^2-4z^3+4z^4+6z^5+2z^6$ so that
$$
{1\over P_{HKR}(z)}=(1+1/z)/HKR^!(z)-(1-5z^2-4z^3+4z^4+6z^5+2z^6)/z
$$
so that
$$
P_{HKR}(z)={1-z \over (1+z)^2(1-2z-z^2)(1-z-z^2+z^3+z^4)}
$$
But
$$
P_R(z)={(1+z)^3(1-z) \over (1-2z-z^2)(1-z-z^2+z^3+z^4)}
$$
so the Avramov spectral sequence degenerates.
Furthermore ${\rm Ext}^*_{HKR}(k,k)$ is not finitely generated.
\medskip
More details: {\bf Case III}
Now the ring is $R =k[x,y,z,u,v]/(xz+yu,zu+uv,yv,z^2,v^2)$. We have that
$$
R^!(z)=\left({1+z\over 1-z-z^2}\right)^2\prod_{n=1}^\infty {1+z^{2n-1}\over 1-z^{2n}}
$$
and
$$
R(z)={1+3z+z^2-2z^3+z^4 \over (1-z)^2}
$$
and the formula (it is still valid)
$$
{1 \over P_R(z)}=(1+1/z)/R^!(z)-R(-z)/z
$$
gives $P_R(z)$.

On the other hand {\tt DGAlgebras} gives that $HKR$ has five generators $X_1,X_2,X_3,X_4,X_5$ of degree
(1,2), three generators $X_5,X_7,X_8$ of degree (2,3), one generator $X_9$ of degree (4,6) and one
 generator $X_{10}$ of degree (5,7). Furthermore, these generators satisfy 29 quadratic relations
of which only the following are non-monomial relations:
$$
X_1X_3+X_3X_4,X_1X_2+X_3X_5,X_1X_8-X_4X_8,X_1X_6+X_3X_7+X_5X_8
$$
Since $IHKR^3=0$ and 
$$
IHKR^2=(X_3X_5,X_2X_5,X_1X_5,X_3X_4,X_2X_4,X_1X_4,X_2X_3,X_5X_8,X_4X_8,X_2X_8,X_3X_7,
$$
$$
X_2X_7,X_1X_7,X_5X_6,X_4X_6,X_3X_6,X_2X_6,X_7X_8,X_6X_8,X_6X_7,X_6X_6)
$$

 It follows that $HKR^!$ is the quotient of the free algebra on 
the ten generators $Y_i$ (which are dual to the $X_i$) by the twosided ideal generated
by the 21 quadratic relations
$$
[Y_1,Y_4],[Y_1,Y_5],[Y_2,Y_3],[Y_2,Y_4],[Y_2,Y_5],[Y_1,Y_2]-[Y_3,Y_5],[Y_1,Y_3]-[Y_3,Y_4],
$$
$$
[Y_1,Y_7],[Y_2,Y_6],[Y_2,Y_7],[Y_2,Y_8],[Y_3,Y_6],[Y_4,Y_6],[Y_5,Y_6],[Y_1,Y_6]-[Y_3,Y_7],
$$
$$
[Y_3,Y_7]-[Y_5,Y_8],[Y_1,Y_8]-[Y_5,Y_8],Y_6Y_6,[Y_6,Y_7],[Y_6,Y_8],[Y_7,Y_8]
$$
where $[,]$ again means the graded commutator so that $[Y_i,Y_j]=Y_iY_j-Y_jY_i$ if
$1 \leq i \leq 5$ and $1 \leq j \leq 5$, Furthermore $[Y_i,Y_j]=Y_iY_j+Y_jY_i$ otherwise.
Note that here the generators should have the degree 2,3,4,5 and we again apply the Ufnarovski trick
introducing a new extra variable T and replacing everywhere $Y_i$ by $T^{deg Y_i-1}Z_i$ gives a new algebra
$A$ with generators of degree 1, such that 
$$
1/HKR^!(z) = 1/A(z)+11z-5z^2-3z^3-x^5-z^6
$$
Again {\tt BERGMAN} give $A(z)$ easily up to degree 20 and we get:
$$
1/HKR^!(z)=1-5z^2-3z^3+7z^4+9z^5+3z^6-5z^7-13z^8-11z^9+2z^{10}+9z^{11}+z^{12}-5z^{13}+
$$
$$
+2z^{14}+7z^{15}+4z^{16}+5z^{17}+6z^{18}-2*z^{19}-8z^{20} \ldots
$$
and since $HKR(x,y)=1+5xy+3xy^2+xy^4+xy^5+7x^2y^2+10x^2y^3+4x^2y^4$
we get $HKR(-z,z)=1-5z^2-3z^3+7z^4+9z^5+3z^5$
and the L\"ofwall formula
$$
1/\Phi_R(z,z)=1/P_{HKR}(z,z) =(1+1/z)/HKR^!(z)-HKR(-z,z)/z
$$
gives
$$
\Phi_R(z,z)=1+5z^2+3z^3+18z^4+21z^5+66z^6+111z^7+274z^8+549z^9+1251*z^{10} \ldots
$$
and the Avramov spectral sequence degenerates up to degree 20 (and probably up to degree $\infty$).

\mysec{7. References}
\medskip
{\frenchspacing\raggedbottom
\def\myref#1{\smallskip\item{[{\bf #1}]}}

\myref{A} D.J.  Anick, {\it A counterexample to a conjecture of Serre},
Ann. of Math. {\bf 115}, 1982, pp. 1-33; {\it Correction}, Ann. of Math. {\bf 116} ,
p. 661.

\myref{A-G} Avramov, L. L.; Golod, E. S.,
{\it The homology of algebra of the Koszul complex of a local Gorenstein ring. (Russian)}, 
Mat. Zametki {\bf 9}, (1971), pp. 53-58.  English translation [Math. Notes {\bf 9} (1971), 30-32].

\myref {AV1} Avramov, L. L.,{\it Obstructions to the existence of multiplicative structures on minimal free resolutions.}
 Amer. J. Math., {\bf 103} (1981), pp. 1-31. 

\myref {AV2} Avramov, L. L.,
{\it Infinite free resolutions,} in
 Six lectures on commutative algebra (Bellaterra, 1996), 
Progr. Math., 166, Birkhäuser, Basel, 1998, pp. 1-118.

\myref {BA} Backelin, J. et al, {\tt BERGMAN},{\it A programme for (non-commutative) Gr\"obner basis calculations} available at {\tt /http://servus.math.su.se/bergman/}
 
\myref{Ei-K} D. Eisenbud, J. Koh,
{\it Nets of alternating matrices and the linear syzygy conjectures},
Adv. Math.{\bf 106} (1994), pp. 1-35.

\myref{GOV} Govorov, V. E. {\it Graded algebras. (Russian)}
Mat. Zametki {\bf 12}, (1972), pp. 197-204.  English translation [Math. Notes {\bf 12} (1973), 556-562].

\myref{KAT} L. Katth\"an, {\it A non-Golod ring with a trivial product on its Koszul homology},

{\tt http://arxiv.org/pdf/1511.04883.pdf}

\myref{LE} J.-M. Lemaire, {\it Alg\`ebres Connexes et Homologie des Espaces
de Lacets}, Lecture Notes in Mathematics, vol. {\bf 422}, Springer-Verlag,
1974.

\myref{L1} L\"ofwall, C., {\it On the subalgebra generated by
one-dimensional elements in the Yoneda {\rm Ext}--algebra}, in
Algebra, algebraic topology, and their interactions,
(       J.--E. Roos, ed),       Lecture Notes in Math., vol.{\bf  1183}, Springer-Verlag,
Berlin--New York, 1986, pp. 291-338.
\myref{L2} C. L\"ofwall, {\it The Yoneda {\rm Ext}--algebra for an equi-characteristic local ring
$(R,m)$ with $m^3=0$}, unpublished manuscript, circa 1976.
\myref{LR1} C. L\"ofwall and J.-E. Roos, {\it Cohomologie des alg\`ebres de Lie
gradu\'ees et s\'eries de Poincar\'e-Betti non rationnelles},
 C.R. Acad. Sc. Paris {\bf 290}, 1980, pp. A733-A736.
\myref{LR2} C. L\"ofwall and J.-E. Roos, {\it  A nonnilpotent 1-2-presented graded Hopf algebra whose Hilbert series converges in the unit circle},
 Adv. Math.{\bf  130}, 1997, pp. 161-200.

\myref{MOO} Moore, Frank, {\it DGAlgebras.} A package for {\tt Macaulay2} that is used to define and manipulate DG algebras,
available at 
{\tt http://www.math.uiuc.edu/Macaulay2/Packages/}

\myref{R1} J.-E. Roos, {\it Relations between the Poincar\'e-Betti series
of Loop Spaces and of Local rings},  Springer Lecture Notes in Math.{\bf 740},
1979, 285-322.
\myref{R2} J.-E. Roos, {\it On the use of graded Lie algebras in the theory of
local rings}, Commutative algebra: Durham 1981 (R. Y. Sharp, ed.) London Math.
Soc. Lecture Notes Ser. vol. 72, Cambridge Univ. Press, Cambridge, 1982, pp. 204–230.
\myref{R3} J.-E. Roos, {\it A computer-aided study of the graded Lie-algebra
of a local commutative noetherian ring}
 (with an Appendix
by Clas L\"ofwall),  Journal of Pure and
Applied Algebra {\bf 91 }, 1994, pp. 255-315.
\myref{U} Ufnarovskij, V., {\it Combinatorial and Asymptotic Methods in Algebra,}
 
in Encyclopaedia of Mathematical Sciences, vol. 57 , Algebra VI

(A.I. Kostrikin and I.R. Shafarevich, eds.), Springer, Berlin, 1994, pp. 1-196.
}

\bigskip
\bigskip
\bigskip
\vfill \break 

%
%
\font\bold=cmbx10 at 14pt

\centerline{\bold APPENDIX 1 }
\medskip
\centerline{\bold A Description of the Homological Behaviour of }
\medskip 
\centerline{\bold  families of quadratic forms in four variables.}
\medskip

\centerline{Jan-Erik Roos}
\centerline{Department of Mathematics}
\centerline{University of Stockholm}
\centerline{S--106 91 Stockholm, SWEDEN}
\centerline{{\tt email: jeroos@math.su.se}}
\bigskip

\def\mysec#1{\bigskip\centerline{\bf #1}\nobreak\par}

\def\cite#1{~[{\bf #1}]}

\bigskip
This is Appendix 1 to our present paper (2015)
about the homological properties of the homology of the Koszul algebra of $codim \geq 4$ local rings.
Part of this appendix was originally included as pages 86-95 in the volume ``Syzygies and Geometry'', october 7-8, 1995, AMS Special Session and
International Conference, Northeastern University, Boston 1995, edited by Antony Iarrobino, Alex Martsinkovsky and Jerzy Weyman.
But that volume had limited circulation, and here we present an update where furthermore
hopefully all misprints in the tables have been corrected (I thank Aldo Conca for help with this).
\bigskip
Let $k$ be a field (for simplicity we assume here that it is of characteristic $0$),
$k[X_1,\ldots,X_n]$ the commutative polynomial ring in $n$ variables over $k$ and let
$f_1,\ldots,f_t$ be quadratic forms in the $X_i$:s. Consider the quotient ring
$$
 R=k[X_1,\ldots,X_n]/(f_1,\ldots f_r) \leqno(1)
$$
which is a graded vector space over $k$.
If $V$ is a vector space over $k$, we denote by $|V|$ the dimension dim$_k(V)$ of $V$ over $k$.
Let $R(t)=\sum_{i=0}^{\infty}|R^i|t^i$ be the Hilbert series of $R$.
In Tables 1-2 below you find a description of the $60$ possible series $R(t)$ when $n=4$
(the ``example'' column in Tables 1-2 refers to the examples in Tables 3-7 below).

Let ${\rm Ext}^*_R(k,k) = {\rm Ext}^{*,*}_R(k,k)$ be the Yoneda Ext-algebra of $R$.
This space is dual to ${\rm Tor}_{*,*}^R(k,k)$ and the second grading comes from
the grading of $R$. In Tables 3-7 you find the possible $83$ series
$$
P_R(x,y) = \sum_{i,j\geq 0}|{\rm Tor}^R_{i,j}(k,k)|x^iy^j \leqno(2)
$$  
for $n=4$. Indeed, these Tables 3-7 should be read in conjunction with Tables 1-2 and Table 8
which describes the relevant non Koszul cases. Let me be more precise:
The sub-algebra of  ${\rm Ext}^*_R(k,k)$ generated by  ${\rm Ext}^1_R(k,k)$ is denoted by $R^!$
and its Hilbert series by $R^!(t)$. For so-called Koszul algebras $R$ we have
$P_R(x,y)=R^!(xy)$ which implies in particular (put $x=-1,y=t$) that
$$
R(t)R^!(-t) = 1. \leqno(3)
$$  
In all cases for $n\leq 4$, {\it except one case}, we have the formula
$$
P_R(x,y)^{-1} = (1+1/x)/R^!(xy)-R(-xy)/x, \leqno(4)
$$
and in the only exceptional case (case {\bf 12} of Table 3) we have
$$
P_R(x,y)^{-1} = (1-1/x^2)/R^!(xy)+R(-xy)/x^2. \leqno(5)
$$
The Tables 3-7 give the Betti numbers as the programme MACAULAY presents them;
thus the first horizontal line for each case gives the $|{\rm Tor}_{i,i}^R(k,k)|$,
the second line gives the $|{\rm Tor}_{i,i+1}^R(k,k)|$, etc. If the ring is a
Koszul algebra there is only one horizontal line and if $R(t)$ is given as
e.g. $H_2$ (which happens in case {\bf 81}), then the corresponding $R^!(t)$
(which is determined by (3)) is denoted by e.g. $A_2$. In the non Koszul cases
we have special notations for the $R^!(t)$ which are given explicitly in Table 8,
and we use the formula (4) and (in case {\bf 12}) the formula (5) to get $P_R(x,y)$.

Among these $83$ cases there are $68$ depth $0$ cases which are new (they are given
in {\bf boldface} in Tables 2-7) and $15$ cases of positive depth which correspond to
$n\leq 3$ and which were known before [1].

Furthermore ${\rm Ext}^{*,*}_R(k,k)$ and its subalgebra $R^!$ are not only Hopf algebras,
but also enveloping algebras of graded Lie algebras $g^{*,*}$ and $\eta^*$. In
[3,4] the relations between the $g^{*,*}$ and $\eta^*$ are explained.
The Lie algebras $\eta^*$ that correspond to Table 8 are either nilpotent 
(maximal degree of nilpotency is $4$, which means that $\eta^5=0$ but  $\eta^4 \neq 0$)
or nice extensions of free Lie algebras (or product of two free Lie algebras) by
nilpotent Lie algebras.

{\it Remarks about the proofs:}
In all the $83$ cases the corresponding dual algebra $R^!$ can be written down
and it is a quotient of a free algebra in $4$ variables with an ideal generated
by quadratic elements. We now want to determine the Hilbert series of that quotient, i.e. $R^!(t)$.
For this we use the programme BERGMAN which calculates a Gr{\"o}bner basis of that ideal.
It turns out that in all the 83 cases we can find a suitable permutation of the variables
for which this ideal has a {\it finite} Gr{\"o}bner basis. Then it follows that the
quotient has a rational series (Govorov and Backelin). Let me illustrate this
explicitly for the most complicated case in table 8, namely case $37$.
Here $R_{37}=k[x,y,z,u]/(x^2,y^2,xy-zu,yz-xu,(x-y)(z-u))$.
It follows that the Koszul dual $R_{37}^!$ is given by the formula:
 $$
R_{37}^! = {k<X,Y,Z,U>\over (Z^2,U^2,XY+YX+ZU+UZ,XZ+ZX-YU-UY,2XZ+2ZX+YZ+ZY+XU+UX)}
$$
where $k<X,Y,Z,U>$ is the free associative algebra in the dual variables
$X,Y,Z,U$.
When we plug this into BERGMAN we see that the correspondig ideal
has an {\it infinite} Gr{\"o}bner basis. However, using the {\tt permutebergman}
addition by Backelin to BERGMAN as described in [5] we see that the order of the variables
(X,U,Z,Y) gives a {\it finite} Gr\"obner basis with 5 elements in degree 2,
2 elements in degree 3, 2 in degree 4, 1 in degree 5 and 1 in degree 6.
Therefore the Hilbert series is rational and can be explicitly determined.
Furthermore it is easy to see that the conditions $L_3$ are satisfied in all cases,
except case 12 of Table 3, where we have $L_4$. This explains the formulae (4) and (5).
Our results are also related to the rational homotopy Lie algebras of CW-complexes [2],[6].

\bigskip

\centerline {BIBLOGRAPHY}

[1] Backelin,J. and Fr\"oberg,R., {\it Poincar\'e series of short artinian rings}, J. of Algebra,
{\bf 96}, 1985, pp. 495-498.

[2] J.-M. Lemaire, {\it Alg\`ebres Connexes et Homologie des Espaces
de Lacets}, Lecture Notes in Mathematics, vol. {\bf 422}, Springer-Verlag,
1974.

[3] L\"ofwall,C. {\it On the subalgebra generated by
one-dimensional elements in the Yoneda {\rm Ext}--algebra}
Algebra, algebraic topology, and their interactions,
(       J.--E. Roos, ed),       Lecture Notes in Math., vol. 1183, Springer-Verlag,
Berlin--New York, 1986, pp. 291-338.

[4] J.-E. Roos, {\it A computer-aided study of the graded Lie-algebra
of a local commutative noetherian ring} (with an Appendix
by Clas L\"ofwall),
Journal of Pure and
Applied Algebra {\bf 91},1994 pp. 255--315.

[5] J.-E. Roos, {\it Three-dimensional manifolds, skew-Gorenstein rings and their cohomology.} J. Commut. Algebra {\bf }2 (2010), no. 4, pp 473---499. 

[6] J.-E. Roos, {\it Relations between the Poincar\'e-Betti series
of Loop Spaces and of Local rings},  Springer Lecture Notes in Math.{\bf 740},
1979, 285-322.

\newdimen\tempdim                
\newdimen\othick \othick=.4pt    
\newdimen\ithick \ithick=.4pt    
\newdimen\spacing \spacing=9pt   
\newdimen\abovehr \abovehr=6pt   
\newdimen\belowhr \belowhr=8pt   
\newdimen\nexttovr \nexttovr=8pt 

\def\r{\hfil&\omit\vrsp\vrule width\othick\cr&}   
\def\rr{\hfil\down{\abovehr}&\omit\vrsp\vrule width\othick\cr
 \noalign{\hrule height\ithick}\up{\belowhr}&}     
\def\up#1{\tempdim=#1\advance\tempdim by1ex
  \vrule height\tempdim width0pt depth0pt} 
\def\down#1{\vrule height0pt depth#1 width0pt} 
\def\large#1#2{\setbox0=\vtop{\hsize#1 \lineskiplimit=0pt \lineskip=1pt
  \baselineskip\spacing \advance\baselineskip by 3pt \noindent
  #2}\tempdim=\dp0\advance\tempdim by\abovehr\box0\down{\tempdim}}
\def\vrsp{\hskip\nexttovr\relax} 
\def\toprule#1{\def\startrule{\hrule height#1\relax}} 
\toprule{\othick}     
\def\nstrut{\vrule height\spacing depth3.5pt width0pt}
\def\exclaim{\char`\!}   
\def\preamble#1{\def\startup{#1}}  
\preamble{&##}  
{\catcode`\!=\active
\gdef!{\hfil\vrule width0pt\vrsp\vrule width\ithick\relax\vrsp&}}

\def\table #1{\vbox\bgroup \setbox0=\hbox{#1}
 \vbox\bgroup\offinterlineskip \catcode`\!=\active
\halign\bgroup##\vrule width\othick\vrsp&\span\startup\nstrut\cr
\noalign{\medskip}
\noalign{\startrule}\up{\belowhr}&}

\def\caption #1{\down{\abovehr}&\omit\vrsp\vrule width\othick\cr
\noalign{\hrule height\othick}\egroup\egroup \setbox1=\lastbox
\tempdim=\wd1 \hbox to\tempdim{\hfil \box0 \hfil} \box1 \smallskip
\hbox to\tempdim{\advance\tempdim by-20pt\hfil\vbox{\hsize\tempdim
\noindent #1}\hfil}\egroup}


%
%
\font\bold=cmbx10 at 14pt

\newdimen\tempdim                
\newdimen\othick \othick=.4pt    
\newdimen\ithick \ithick=.4pt    
\newdimen\spacing \spacing=9pt   
\newdimen\abovehr \abovehr=6pt   
\newdimen\belowhr \belowhr=8pt   
\newdimen\nexttovr \nexttovr=8pt 

\def\r{\hfil&\omit\vrsp\vrule width\othick\cr&}   
\def\rr{\hfil\down{\abovehr}&\omit\vrsp\vrule width\othick\cr
 \noalign{\hrule height\ithick}\up{\belowhr}&}     
\def\up#1{\tempdim=#1\advance\tempdim by1ex
  \vrule height\tempdim width0pt depth0pt} 
\def\down#1{\vrule height0pt depth#1 width0pt} 
\def\large#1#2{\setbox0=\vtop{\hsize#1 \lineskiplimit=0pt \lineskip=1pt
  \baselineskip\spacing \advance\baselineskip by 3pt \noindent
  #2}\tempdim=\dp0\advance\tempdim by\abovehr\box0\down{\tempdim}}
\def\vrsp{\hskip\nexttovr\relax} 
\def\toprule#1{\def\startrule{\hrule height#1\relax}} 
\toprule{\othick}     
\def\nstrut{\vrule height\spacing depth3.5pt width0pt}
\def\exclaim{\char`\!}   
\def\preamble#1{\def\startup{#1}}  
\preamble{&##}  
{\catcode`\!=\active
\gdef!{\hfil\vrule width0pt\vrsp\vrule width\ithick\relax\vrsp&}}

\def\table #1{\vbox\bgroup \setbox0=\hbox{#1}
 \vbox\bgroup\offinterlineskip \catcode`\!=\active
\halign\bgroup##\vrule width\othick\vrsp&\span\startup\nstrut\cr
\noalign{\medskip}
\noalign{\startrule}\up{\belowhr}&}

\def\caption #1{\down{\abovehr}&\omit\vrsp\vrule width\othick\cr
\noalign{\hrule height\othick}\egroup\egroup \setbox1=\lastbox
\tempdim=\wd1 \hbox to\tempdim{\hfil \box0 \hfil} \box1 \smallskip
\hbox to\tempdim{\advance\tempdim by-20pt\hfil\vbox{\hsize\tempdim
\noindent #1}\hfil}\egroup}

$$\table{\bf Table 1. THE FIRST 30 POSSIBLE HILBERT SERIES $H_{1},\ldots,H_{30}$}
   !\hfil Rational form !\multispan{7} \hfil Taylor series (Hilbert function) !\hfil Example \rr
\hfil $H_1$ ! $1+4\,t$ !\hfil 1 ! 4 ! 0 !  0 ! 0 ! 0 ! 0 !  \hfil 83\rr
\hfil $H_2$ ! $1+4\,t+t^2$ !\hfil 1 ! 4 !  1 ! 0 ! 0 ! 0 ! 0 ! \hfil 81  \rr
\hfil $H_3$ ! $(1+3\,t-3\,t^2)/(1-t)$ !\hfil 1 ! 4 !  1 ! 1 ! 1 ! 1 ! 1 !\hfil 82   \rr
\hfil $H_4$ ! $1+4\,t+2\,t^2$ !\hfil 1 ! 4 ! 2 !  0 ! 0 ! 0 ! 0 ! \hfil 78 \rr
\hfil $H_5$ ! $(1+3\,t-2\,t^2-t^3)/(1-t)$ !\hfil 1 ! 4 !  2 ! 1 ! 1 ! 1 ! 1 !\hfil 79   \rr
\hfil $H_6$ ! $(1+3\,t-2\,t^2)/(1-t)$ !\hfil 1 ! 4 !  2 ! 2 ! 2 ! 2 ! 2 ! \hfil 80  \rr
\hfil $H_7$ ! $1+4\,t+3\,t^2=(1+t)(1+3\,t)$ !\hfil 1 ! 4 ! 3 !  0 ! 0 ! 0 ! 0 ! \hfil 71 \rr
\hfil $H_8$ ! $1+4\,t+3\,t^2+t^3$ !\hfil 1 ! 4 ! 3 !  1 ! 0 ! 0 ! 0 ! \hfil 73 \rr
\hfil $H_9$ ! $(1+3\,t-t^2-2\,t^3)/(1-t)$ !\hfil 1 ! 4 ! 3 !  1 ! 1 ! 1 ! 1 !\hfil 72,74 \rr
\hfil $H_{10}$ ! $(1+3\,t-t^2-t^3)/ (1-t)$ !\hfil 1 ! 4 ! 3 !  2 ! 2 ! 2 ! 2 ! \hfil 75 \rr
\hfil $H_{11}$ ! $(1+3\,t-t^2)/(1-t)$ !\hfil 1 ! 4 ! 3 !  3 ! 3 ! 3 ! 3 ! \hfil 76 \rr
\hfil $H_{12}$ ! $(1+2\,t-4\,t^2+2\,t^3)/(1-t)^2$ !\hfil 1 ! 4 ! 3 !  4 ! 5 ! 6 ! 7 ! \hfil 77 \rr
\hfil $H_{13}$ ! $1+4\,t+4\,t^2=(1+2\,t)^2$ !\hfil 1 ! 4 ! 4 !  0 ! 0 ! 0 ! 0 ! \hfil 54,55,56,57 \rr
\hfil $H_{14}$ ! $\scriptstyle{1+4\,t+4\,t^2+t^3=(1+t)(1+3\,t+t^2)}$ !\hfil 1 ! 4 ! 4 !  1 ! 0 ! 0 ! 0 !\hfil 61  \rr
\hfil $H_{15}$ ! $(1+3\,t-3\,t^3)/(1-t)$ !\hfil 1 ! 4 ! 4 !  1 ! 1 ! 1 ! 1 !\hfil 58,59,60,62  \rr
\hfil $H_{16}$ ! $(1+3\,t-2\,t^3-t^4)/(1-t)$ !\hfil 1 ! 4 ! 4 !  2 ! 1 ! 1 ! 1 ! \hfil 64 \rr
\hfil $H_{17}$ ! $(1+2\,t-2\,t^2)(1+t)/(1-t)$ !\hfil 1 ! 4 ! 4 !  2 ! 2 ! 2 ! 2 ! \hfil 63,65 \rr
\hfil $H_{18}$ ! $(1+3\,t-t^3)/(1-t)$ !\hfil 1 ! 4 ! 4 !  3 ! 3 ! 3 ! 3 ! \hfil 66 \rr
\hfil $H_{19}$ ! $(1+3\,t)/(1-t)$ !\hfil 1 ! 4 ! 4 !  4 ! 4 ! 4 ! 4 !  \hfil 68\rr
\hfil $H_{20}$ ! $(1+2\,t-3\,t^2+t^4)/(1-t)^2$ !\hfil 1 ! 4 ! 4 !  4 ! 5 ! 6 ! 7 ! \hfil 67,69 \rr
\hfil $H_{21}$ ! $(1+2\,t-3\,t^2+t^3)/(1-t)^2$ !\hfil 1 ! 4 ! 4 !  5 ! 6 ! 7 ! 8 ! \hfil 70 \rr
\hfil $H_{22}$ ! $1+4\,t+5\,t^2$ !\hfil 1 ! 4 ! 5 !  0 ! 0 ! 0 ! 0 !  \hfil 29 \rr
\hfil $H_{23}$ ! $1+4\,t+5\,t^2+t^3$ !\hfil 1 ! 4 ! 5 !  1 ! 0 ! 0 ! 0 ! \hfil 30,32 \rr
\hfil $H_{24}$ ! $(1+3\,t+t^2-4\,t^3)/(1-t)$ !\hfil 1 ! 4 ! 5 !  1 ! 1 ! 1 ! 1 !\hfil 31,33 \rr
\hfil $H_{25}$ ! $\scriptstyle{1+4\,t+5\,t^2+2\,t^3=(1+t)^2(1+2\,t)}$ !\hfil 1 ! 4 ! 5 !  2 ! 0 ! 0 ! 0 ! \hfil 38 \rr
\hfil $H_{26}$ ! $(1+3\,t+t^2-3\,t^3-t^4)/(1-t)$ !\hfil 1 ! 4 ! 5 !  2 ! 1 ! 1 ! 1 !\hfil 39  \rr
\hfil $H_{27}$ ! $(1+3\,t+t^2-3\,t^3)/(1-t)$ !\hfil 1 ! 4 ! 5 !  2 ! 2 ! 2 ! 2 ! \hfil 34,35,36,37,40 \rr
\hfil $H_{28}$ ! $(1+2\,t-t^2-t^3)(1+t)/(1-t)$ !\hfil 1 ! 4 ! 5 !  3 ! 2 ! 2 ! 2 ! \hfil 43 \rr
\hfil $H_{29}$ ! $(1+2\,t)(1+t-t^2)/(1-t)$ !\hfil 1 ! 4 ! 5 !  3 ! 3 ! 3 ! 3 !\hfil 41,42,44 \rr
\hfil $H_{30}$ ! $(1+2\,t-t^2)(1+t)/(1-t)$ !\hfil 1 ! 4 ! 5 !  4 ! 4 ! 4 ! 4 ! \hfil 46
\caption{\sl On the next page you can see the last 30 possible Hilbert series.} 
$$
%
\font\bold=cmbx10 at 14pt

\newdimen\tempdim                
\newdimen\othick \othick=.4pt    
\newdimen\ithick \ithick=.4pt    
\newdimen\spacing \spacing=9pt   
\newdimen\abovehr \abovehr=6pt   
\newdimen\belowhr \belowhr=8pt   
\newdimen\nexttovr \nexttovr=8pt 

\def\r{\hfil&\omit\vrsp\vrule width\othick\cr&}   
\def\rr{\hfil\down{\abovehr}&\omit\vrsp\vrule width\othick\cr
 \noalign{\hrule height\ithick}\up{\belowhr}&}     
\def\up#1{\tempdim=#1\advance\tempdim by1ex
  \vrule height\tempdim width0pt depth0pt} 
\def\down#1{\vrule height0pt depth#1 width0pt} 
\def\large#1#2{\setbox0=\vtop{\hsize#1 \lineskiplimit=0pt \lineskip=1pt
  \baselineskip\spacing \advance\baselineskip by 3pt \noindent
  #2}\tempdim=\dp0\advance\tempdim by\abovehr\box0\down{\tempdim}}
\def\vrsp{\hskip\nexttovr\relax} 
\def\toprule#1{\def\startrule{\hrule height#1\relax}} 
\toprule{\othick}     
\def\nstrut{\vrule height\spacing depth3.5pt width0pt}
\def\exclaim{\char`\!}   
\def\preamble#1{\def\startup{#1}}  
\preamble{&##}  
{\catcode`\!=\active
\gdef!{\hfil\vrule width0pt\vrsp\vrule width\ithick\relax\vrsp&}}

\def\table #1{\vbox\bgroup \setbox0=\hbox{#1}
 \vbox\bgroup\offinterlineskip \catcode`\!=\active
\halign\bgroup##\vrule width\othick\vrsp&\span\startup\nstrut\cr
\noalign{\medskip}
\noalign{\startrule}\up{\belowhr}&}

\def\caption #1{\down{\abovehr}&\omit\vrsp\vrule width\othick\cr
\noalign{\hrule height\othick}\egroup\egroup \setbox1=\lastbox
\tempdim=\wd1 \hbox to\tempdim{\hfil \box0 \hfil} \box1 \smallskip
\hbox to\tempdim{\advance\tempdim by-20pt\hfil\vbox{\hsize\tempdim
\noindent #1}\hfil}\egroup}

$$\table{\bf Table 2. THE LAST 30 POSSIBLE HILBERT SERIES $H_{31},\ldots,H_{60}$ }
   !\hfil Rational factored form !\multispan{7} \hfil Taylor series (Hilbert function) !\hfil Example \rr
\hfil $H_{31}$ ! $(1+2\,t-2\,t^2-2\,t^3+2\,t^4)/(1-t)^2$ !\hfil 1 ! 4 ! 5 !  4 ! 5 ! 6 ! 7 ! \hfil 45,47,48 \rr
\hfil $H_{32}$ ! $(1+3\,t+t^2)/(1-t)$ !\hfil 1 ! 4 ! 5 !  5 ! 5 ! 5 ! 5 ! \hfil 50 \rr
\hfil $H_{33}$ ! $(1+2\,t-2\,t^2-t^3+t^4)/(1-t)^2$ !\hfil 1 ! 4 ! 5 !  5 ! 6 ! 7 ! 8 !\hfil 49,51 \rr
\hfil $H_{34}$ ! $(1+2\,t-2\,t^2)/(1-t)^2$ !\hfil 1 ! 4 ! 5 !  6 ! 7 ! 8 ! 9 ! \hfil 52 \rr
\hfil $H_{35}$ ! $(1+2\,t-2\,t^2+t^3)/(1-t)^2$ !\hfil 1 ! 4 ! 5 !  7 ! 9 ! 11 ! 13 ! \hfil 53  \rr
\hfil $H_{36}$ ! $1+4\,t+6\,t^2+4\,t^3+t^4$ !\hfil 1 ! 4 ! 6 !  4 ! 1 ! 0 ! 0 !  \hfil 11\rr 
\hfil $H_{37}$ ! $(1+3\,t+2\,t^2-2\,t^3-3\,t^4)/(1-t)$ !\hfil 1 ! 4 ! 6 !  4 ! 1 ! 1 ! 1 ! \hfil 12 \rr 
\hfil $H_{38}$ ! $(1+3\,t+2\,t^2-2\,t^3-2\,t^4)/(1-t)$ !\hfil 1 ! 4 ! 6 !  4 ! 2 ! 2 ! 2 ! \hfil 13 \rr 
\hfil $H_{39}$ ! $(1+3\,t+2\,t^2-2\,t^3-t^4)/(1-t)$ !\hfil 1 ! 4 ! 6 !  4 ! 3 ! 3 ! 3 !  \hfil 14\rr 
\hfil $H_{40}$ ! $(1+3\,t+2\,t^2-2\,t^3)/(1-t)$ !\hfil 1 ! 4 ! 6 !  4 ! 4 ! 4 ! 4 !  \hfil 15\rr 
\hfil $H_{41}$ ! $(1+2\,t-t^2-4\,t^3+3\,t^4)/(1-t)^2$ !\hfil 1 ! 4 ! 6 !  4 ! 5 ! 6 ! 7 !  \hfil 16\rr 
\hfil $H_{42}$ ! $(1+t-t^2)(1+t)^2/(1-t)$ !\hfil 1 ! 4 ! 6 !  5 ! 4 ! 4 ! 4 !  \hfil 18\rr 
\hfil $H_{43}$ ! $(1+3\,t+2\,t^2-t^3)/(1-t)$ !\hfil 1 ! 4 ! 6 !  5 ! 5 ! 5 ! 5 !  \hfil 17,19\rr 
\hfil $H_{44}$ ! $(1+2\,t-t^2-3\,t^3+2\,t^4)/(1-t)^2$ !\hfil 1 ! 4 ! 6 !  5 ! 6 ! 7 ! 8 !  \hfil 20\rr
\hfil $H_{45}$ ! $(1+t)(1+2\,t)/(1-t)$ !\hfil 1 ! 4 ! 6 !  6 ! 6 ! 6 ! 6 !  \hfil 23\rr
\hfil $H_{46}$ ! $(1+t-t^2)^2/(1-t)^2$ !\hfil 1 ! 4 ! 6 !  6 ! 7 ! 8 ! 9 !  \hfil 21,22,24\rr
\hfil $H_{47}$ ! $(1+2\,t-t^2-t^3)/(1-t)^2$ !\hfil 1 ! 4 ! 6 !  7 ! 8 ! 9 ! 10 !  \hfil 26\rr
\hfil $H_{48}$ ! $(1+t-2\,t^2+t^3)(1+t)/(1-t)^2$ !\hfil 1 ! 4 ! 6 ! 7 ! 9 ! 11 ! 13 ! \hfil 25 \rr
\hfil $H_{49}$ ! $(1+2\,t-t^2)/(1-t)^2$ !\hfil 1 ! 4 ! 6 !  8 ! 10 ! 12 ! 14 !  \hfil 27\rr
\hfil $H_{50}$ ! $(1+t-3\,t^2+3\,t^3-t^4)/(1-t)^3$ !\hfil 1 ! 4 ! 6 !  10 ! 15 ! 21 ! 28 ! \hfil 28 \rr
\hfil $H_{51}$ ! $(1+t)^3/(1-t)$ !\hfil 1 ! 4 ! 7 !  8 ! 8 ! 8 ! 8 !  \hfil 5\rr
\hfil $H_{52}$ ! $(1+2\,t-2\,t^3)/(1-t)^2$ !\hfil 1 ! 4 ! 7 !  8 ! 9 ! 10 ! 11 !  \hfil 6\rr
\hfil $H_{53}$ ! $(1+2\,t-2\,t^3+t^4)/(1-t)^2$ !\hfil 1 ! 4 ! 7 !  8 ! 10 ! 12 ! 14 !  \hfil 7\rr
\hfil $H_{54}$ ! $(1+t-t^2)(1+t)/(1-t)^2$ !\hfil 1 ! 4 ! 7 !  9 ! 11 ! 13 ! 15 !  \hfil 8\rr
\hfil $H_{55}$ ! $(1+2\,t)/(1-t)^2$ !\hfil 1 ! 4 ! 7 !  10 ! 13 ! 16 ! 19 !  \hfil 9\rr
\hfil $H_{56}$ ! $(1+t-2\,t^2+t^3)/(1-t)^3$ !\hfil 1 ! 4 ! 7 !  11 ! 16 ! 22 ! 29 ! \hfil 10 \rr
\hfil $H_{57}$ ! $(1+t)^2/(1-t)^2$ !\hfil 1 ! 4 ! 8 ! 12 ! 16 ! 20 ! 24 ! \hfil 3 \rr
\hfil $H_{58}$ ! $(1+t-t^2)/(1-t)^3$ !\hfil 1 ! 4 ! 8 !  13 ! 19 ! 26 ! 34 ! \hfil 4 \rr
\hfil $H_{59}$ ! $(1+t)/(1-t)^3$ !\hfil 1 ! 4 ! 9 ! 16 ! 25 ! 36 ! 49 !  \hfil 2\rr
\hfil $H_{60}$ ! $1/(1-t)^4$ !\hfil 1 ! 4 ! 10 !  20 ! 35 ! 56 ! 84 !  \hfil 1
\caption{\sl } 
$$ 

%
\font\bold=cmbx10 at 14pt

\newdimen\tempdim                
\newdimen\othick \othick=.4pt    
\newdimen\ithick \ithick=.4pt    
\newdimen\spacing \spacing=9pt   
\newdimen\abovehr \abovehr=6pt   
\newdimen\belowhr \belowhr=8pt   
\newdimen\nexttovr \nexttovr=8pt 

\def\r{\hfil&\omit\vrsp\vrule width\othick\cr&}   
\def\rr{\hfil\down{\abovehr}&\omit\vrsp\vrule width\othick\cr
 \noalign{\hrule height\ithick}\up{\belowhr}&}     
\def\up#1{\tempdim=#1\advance\tempdim by1ex
  \vrule height\tempdim width0pt depth0pt} 
\def\down#1{\vrule height0pt depth#1 width0pt} 
\def\large#1#2{\setbox0=\vtop{\hsize#1 \lineskiplimit=0pt \lineskip=1pt
  \baselineskip\spacing \advance\baselineskip by 3pt \noindent
  #2}\tempdim=\dp0\advance\tempdim by\abovehr\box0\down{\tempdim}}
\def\vrsp{\hskip\nexttovr\relax} 
\def\toprule#1{\def\startrule{\hrule height#1\relax}} 
\toprule{\othick}     
\def\nstrut{\vrule height\spacing depth3.5pt width0pt}
\def\exclaim{\char`\!}   
\def\preamble#1{\def\startup{#1}}  
\preamble{&##}  
{\catcode`\!=\active
\gdef!{\hfil\vrule width0pt\vrsp\vrule width\ithick\relax\vrsp&}}

\def\table #1{\vbox\bgroup \setbox0=\hbox{#1}
 \vbox\bgroup\offinterlineskip \catcode`\!=\active
\halign\bgroup##\vrule width\othick\vrsp&\span\startup\nstrut\cr
\noalign{\medskip}
\noalign{\startrule}\up{\belowhr}&}

\def\caption #1{\down{\abovehr}&\omit\vrsp\vrule width\othick\cr
\noalign{\hrule height\othick}\egroup\egroup \setbox1=\lastbox
\tempdim=\wd1 \hbox to\tempdim{\hfil \box0 \hfil} \box1 \smallskip
\hbox to\tempdim{\advance\tempdim by-20pt\hfil\vbox{\hsize\tempdim
\noindent #1}\hfil}\egroup}
$$\table{Table 3. The first 19 cases (case numbers in {\bf boldface} indicate depth 0 rings). }
Case !\multispan{8} The graded Betti Numbers! ideal gen:s (ex.)  ! $R(t)$ ! $R^{\exclaim}(t)$     \rr
\hfill $1$ ! \hfill 1 !\hfill 4     !\hfill 6   !\hfill 4   !\hfill 1   !\hfill 0     !\hfill 0  !\hfill 0   ! $0$  !$H_{60}$      !$A_{60}$  \rr
\hfill 2 !\hfill 1  !\hfill 4    !\hfill 7  !\hfill 8   !\hfill 8   !\hfill 8    !\hfill 8  !\hfill 8   ! $x^2$  !$H_{59}$     !$A_{59}$  \rr
 \hfill 3   ! \hfill 1    !\hfill    4    !\hfill    8    !\hfill    12   !\hfill     16   !\hfill    20   !\hfill     24   !\hfill     28       !$x^2,y^2$       !$H_{57}$  !$A_{57}$\rr   
\hfill 4 ! \hfill 1  !  \hfill 4    !    \hfill 8   !\hfill   13   !\hfill  21   !\hfill 34     !\hfill 55   !\hfill    89  ! $x^2,xy$   !$H_{58}$   !$A_{58}$ \rr
\hfill 5     !\hfill    1    !\hfill    4    !\hfill    9    !\hfill    16   !\hfill     25    !\hfill    36   !\hfill     49   !\hfill     64      !$x^2,y^2,z^2$       !$H_{51}$!$A_{51}$\rr    
\hfill 6 ! \hfill 1  !  \hfill 4    !    \hfill 9   !\hfill   16   !\hfill  25   !\hfill 36     !\hfill 49   ! \hfill 64    ! $x^2,y^2+xz,yz$    !$H_{52}$  !$A_{51}$  \r 
    ! \hfill -  !\hfill  -    ! \hfill   -   !\hfill    1   !\hfill   6   !\hfill 21     ! \hfill 56  ! \hfill 126      !    !   !  \r 
    ! \hfill -  !  \hfill -    !    \hfill  -  !    \hfill -   !   \hfill -   ! \hfill -     ! \hfill 1  ! \hfill 8      ! !      !  \rr
\hfill {\bf 7} ! \hfill 1  !  \hfill 4    !    \hfill 9   !\hfill   16   !\hfill  25   !\hfill 36     !\hfill 49   !\hfill 64     !$x^2+y^2,z^2+u^2,$   !$H_{53}$    !$A_{51}$  \r 
\hfill  ! \hfill -  !  \hfill -    !    \hfill -   !\hfill   2   !\hfill  12   !\hfill 42     !\hfill 112   ! \hfill 252     !$xz+yu$  !     !  \r 
\hfill  ! \hfill -  !  \hfill -    !    \hfill -   !\hfill   -   !\hfill  -   !\hfill -     !\hfill 4   !  \hfill 32    !     !  !  \rr
\hfill 8     !\hfill    1    !\hfill    4    !\hfill    9    !\hfill    17   !\hfill     30   !\hfill     51   !\hfill     85   !\hfill 140         !$x^2,y^2,xz$      !$H_{54}$   !$A_{54}$\rr    
\hfill 9     !\hfill    1    !\hfill    4    !\hfill    9    !\hfill    18   !\hfill     36   !\hfill     72   !\hfill     144   !\hfill 288         !$x^2,xy,y^2$       !$H_{55}$ !$A_{55}$\rr   
\hfill 10    !\hfill    1    !\hfill    4    !\hfill    9    !\hfill    19   !\hfill     41   !\hfill     88   !\hfill     189   !\hfill 406         !$x^2,xy,xz$      !$H_{56}$  !$A_{56}$\rr  
\hfill {\bf 11}    !\hfill    1    !\hfill    4    !\hfill    10    !\hfill    20    !\hfill    35   !\hfill     56   !\hfill     84   !\hfill 120         !$x^2,y^2,z^2,u^2$      !$H_{36}$   !$A_{36}$\rr  
 
\hfill {\bf 12}    !\hfill    1    !\hfill    4    !\hfill    10   !\hfill     20   !\hfill     35   !\hfill     56   !\hfill     84     !\hfill 120        !$x^2+xy,y^2+xu,$       !$H_{37}$!$A_{36}$\r   
 !\hfill          -    !\hfill    -    !\hfill    -    !\hfill    -    !\hfill    -    !\hfill    -    !\hfill    -    !\hfill    -   !$z^2+xu,u^2+zu$         !        !\r  
     !\hfill          -    !\hfill    -    !\hfill    -    !\hfill    1    !\hfill    7   !\hfill     29    !\hfill    91    !\hfill  239     !\hfill                !\hfill !\r   
     !\hfill           -    !\hfill    -    !\hfill    -    !\hfill    -    !\hfill    -    !\hfill    -    !\hfill    -    !\hfill    -   !\hfill          !\hfill       !\r  
     !\hfill           -    !\hfill    -    !\hfill    -    !\hfill    -    !\hfill    -    !\hfill    -    !\hfill    1    !\hfill  10     !\hfill          !\hfill         !\rr  
\hfill {\bf 13}    !\hfill    1    !\hfill    4   !\hfill     10   !\hfill     20   !\hfill     35   !\hfill     56    !\hfill    84   !\hfill 120        !$x^2+z^2+u^2,$      !$H_{38}$!$A_{36}$\r    
     !\hfill    -   !\hfill    -    !\hfill    -    !\hfill    1    !\hfill    6    !\hfill     22    !\hfill    62   !\hfill   148       !$y^2,xz,yu+zu$       !\hfill !\r    
     !\hfill    -   !\hfill    -    !\hfill    -    !\hfill    -    !\hfill    -    !\hfill    -    !\hfill    1    !\hfill     8        !\hfill       !\hfill !\rr    
\hfill {\bf 14}    !\hfill    1    !\hfill    4   !\hfill     10   !\hfill     20   !\hfill     35   !\hfill     56    !\hfill    84   !\hfill  120        !$xz,y^2,z^2+u^2,$       !$H_{39}$!$A_{36}$\r  
     !\hfill       !\hfill    -    !\hfill    -    !\hfill        2    !\hfill    13   !\hfill     51   !\hfill     153   !\hfill  387        !$yu+zu$       !\hfill !\r    
     !\hfill       !\hfill    -    !\hfill    -    !\hfill       -    !\hfill    -    !\hfill    -     !\hfill    4   !\hfill  36        !\hfill       !\hfill !\rr   
\hfill {\bf 15}    !\hfill    1    !\hfill    4   !\hfill     10   !\hfill     20   !\hfill     35   !\hfill     56    !\hfill    84   !\hfill   120       !$xz,y^2,yz+u^2,$       !$H_{40}$!$A_{36}$\r   
     !\hfill       !\hfill    -    !\hfill    -    !\hfill       3   !\hfill     20   !\hfill     80   !\hfill     244   !\hfill 626         !$yu+zu$      !\hfill !\r    
     !\hfill       !\hfill    -    !\hfill    -    !\hfill       -    !\hfill    -    !\hfill    -     !\hfill    9   !\hfill  84        !\hfill       !\hfill !\rr   
\hfill {\bf 16}    !\hfill    1    !\hfill    4   !\hfill     10    !\hfill    20   !\hfill     35    !\hfill    56    !\hfill    84   !\hfill 120         !$xy+z^2+yu,y^2,$       !$H_{41}$!$A_{36}$\r   
     !\hfill       !\hfill        -    !\hfill    -    !\hfill    4    !\hfill    26   !\hfill     103   !\hfill     312   !\hfill 797         !$yu+zu,xz$       !\hfill !\r     
     !\hfill       !\hfill       -    !\hfill    -    !\hfill    -    !\hfill    -    !\hfill    -     !\hfill    16   !\hfill  144        !\hfill       !\hfill ! \rr  
\hfill {\bf 17}    !\hfill    1    !\hfill    4   !\hfill     10    !\hfill    21    !\hfill    39    !\hfill    66    !\hfill    104    !\hfill    155   !$xz,yz+xu,$          !$H_{43}$!$A_{43.1}$\r   
     !\hfill       !\hfill        -    !\hfill    -    !\hfill    -    !\hfill    1     !\hfill    8    !\hfill    36   !\hfill     121   !$y^2,yu+zu$          !\hfill       ! \rr  
\hfill {\bf 18}    !\hfill    1    !\hfill    4   !\hfill     10    !\hfill    21   !\hfill     40    !\hfill    72    !\hfill    125     !\hfill  212       !$x^2,y^2,z^2,yu$       !$H_{42}$ !$A_{42}$\rr   
\hfill {\bf 19}    !\hfill    1    !\hfill    4   !\hfill     10   !\hfill     21   !\hfill     40    !\hfill    72     !\hfill    125     !\hfill 212        !$xz,y^2,$       !$H_{43}$ !$A_{42}$\r    
     !\hfill       !\hfill    -        !\hfill    -    !\hfill    1    !\hfill    7    !\hfill    29    !\hfill    93     !\hfill   255      !$yu+zu,u^2$       ! !\r   
     !\hfill       !\hfill    -       !\hfill    -    !\hfill    -    !\hfill    -     !\hfill    -    !\hfill    1     !\hfill   10      !\hfill       !\hfill !
\caption{\sl } 
$$
%
\font\bold=cmbx10 at 14pt

\newdimen\tempdim                
\newdimen\othick \othick=.4pt    
\newdimen\ithick \ithick=.4pt    
\newdimen\spacing \spacing=9pt   
\newdimen\abovehr \abovehr=6pt   
\newdimen\belowhr \belowhr=8pt   
\newdimen\nexttovr \nexttovr=8pt 

\def\r{\hfil&\omit\vrsp\vrule width\othick\cr&}   
\def\rr{\hfil\down{\abovehr}&\omit\vrsp\vrule width\othick\cr
 \noalign{\hrule height\ithick}\up{\belowhr}&}     
\def\up#1{\tempdim=#1\advance\tempdim by1ex
  \vrule height\tempdim width0pt depth0pt} 
\def\down#1{\vrule height0pt depth#1 width0pt} 
\def\large#1#2{\setbox0=\vtop{\hsize#1 \lineskiplimit=0pt \lineskip=1pt
  \baselineskip\spacing \advance\baselineskip by 3pt \noindent
  #2}\tempdim=\dp0\advance\tempdim by\abovehr\box0\down{\tempdim}}
\def\vrsp{\hskip\nexttovr\relax} 
\def\toprule#1{\def\startrule{\hrule height#1\relax}} 
\toprule{\othick}     
\def\nstrut{\vrule height\spacing depth3.5pt width0pt}
\def\exclaim{\char`\!}   
\def\preamble#1{\def\startup{#1}}  
\preamble{&##}  
{\catcode`\!=\active
\gdef!{\hfil\vrule width0pt\vrsp\vrule width\ithick\relax\vrsp&}}

\def\table #1{\vbox\bgroup \setbox0=\hbox{#1}
 \vbox\bgroup\offinterlineskip \catcode`\!=\active
\halign\bgroup##\vrule width\othick\vrsp&\span\startup\nstrut\cr
\noalign{\medskip}
\noalign{\startrule}\up{\belowhr}&}

\def\caption #1{\down{\abovehr}&\omit\vrsp\vrule width\othick\cr
\noalign{\hrule height\othick}\egroup\egroup \setbox1=\lastbox
\tempdim=\wd1 \hbox to\tempdim{\hfil \box0 \hfil} \box1 \smallskip
\hbox to\tempdim{\advance\tempdim by-20pt\hfil\vbox{\hsize\tempdim
\noindent #1}\hfil}\egroup}
$$\table{Table 4. Cases 20--37 }
Case !\multispan{8} The graded Betti Numbers! ideal gen:s (ex.)  ! $R(t)$ ! $R^{\exclaim}(t)$     \rr
\hfill {\bf 20}    !\hfill    1    !\hfill    4   !\hfill     10   !\hfill     21   !\hfill     40    !\hfill    72   !\hfill     125     !\hfill 212        !$xz,y^2,$       !$H_{44}$    !$A_{42}$\r  
     !\hfill       !\hfill    -       !\hfill    -    !\hfill    2   !\hfill     13    !\hfill    52    !\hfill    163     !\hfill 440       ! $yu+z^2,$       !  !\r  
     !\hfill       !\hfill    -        !\hfill    -    !\hfill    -    !\hfill    -     !\hfill    -     !\hfill    4     !\hfill 36        ! $yu+zu$       !   !\rr 
\hfill {\bf 21}   !\hfill    1    !\hfill    4    !\hfill    10   !\hfill     22    !\hfill    45    !\hfill    88   !\hfill     167   !\hfill     310     !$xz,y^2,z^2,yu+zu$             !$H_{46}$   !$A_{46}$\rr  
\hfill {\bf 22}   !\hfill    1    !\hfill    4    !\hfill    10   !\hfill     22   !\hfill     45    !\hfill    88   !\hfill     167   !\hfill     311             !$x^2+xy,xu,$      !$H_{46}$ !$A_{46.1}$\r   
       !\hfill    -    !\hfill    -    !\hfill    -    !\hfill    -    !\hfill    -     !\hfill    -     !\hfill    1     !\hfill  7     !$xz+yu,y^2$ !\hfill !\rr   
\hfill  23   !\hfill    1    !\hfill    4   !\hfill     10   !\hfill     22   !\hfill     46    !\hfill    94    !\hfill    190     !\hfill 382        !$xz,xu,y^2,z^2$       !$H_{45}$ !$A_{45}$\rr   
\hfill {\bf 24}    !\hfill    1    !\hfill    4 !\hfill 10   !\hfill     22   !\hfill     46    !\hfill    94    !\hfill    190   !\hfill  382        !$xz,y^2,$       !$H_{46}$ !$A_{45}$\r   
     !\hfill       !\hfill    -        !\hfill    -    !\hfill    1    !\hfill    6    !\hfill    23    !\hfill    72   !\hfill 201         !$yz+z^2,$       !\hfill !\r   
     !\hfill       !\hfill    -        !\hfill    -    !\hfill    -    !\hfill    -     !\hfill    -     !\hfill    1   !\hfill   8       !$yu+zu$       !\hfill !\rr   
\hfill {\bf 25 }   !\hfill    1    !\hfill    4   !\hfill     10   !\hfill     23   !\hfill     51   !\hfill     111   !\hfill     240   !\hfill 517         !$x^2,xy,xz,u^2$       !$H_{48}$ \hfill !$A_{48}$\rr   
\hfill 26    !\hfill    1    !\hfill    4   !\hfill     10   !\hfill     23   !\hfill     52   !\hfill     117   !\hfill     263   !\hfill 591        !$xz,y^2,yu,zu$       !$H_{47}$ !$A_{47}$\rr   
\hfill 27    !\hfill    1    !\hfill    4   !\hfill     10   !\hfill     24   !\hfill     58   !\hfill     140   !\hfill     338   !\hfill 816         !$xy,xz,y^2,yz$       !$H_{49}$  !$A_{49}$\rr  
\hfill {\bf 28}    !\hfill    1    !\hfill    4   !\hfill     10   !\hfill     26   !\hfill     69   !\hfill     181   !\hfill     476   !\hfill 1252         !$x^2,xy,xz,xu$       !$H_{50}$  !$A_{50}$\rr  
\hfill {\bf 29}    !\hfill    1    !\hfill    4   !\hfill     11   !\hfill     24   !\hfill     46    !\hfill    80   !\hfill     130   !\hfill 200         !$x^2+xy,y^2+xu,$ !$H_{22}$   !$A_{22.1}$\r 
       !\hfill    -    !\hfill    -    !\hfill    -    !\hfill    5   !\hfill     36   !\hfill     159   !\hfill     536   !\hfill 1519         !$z^2+xu,zu+u^2,$       !\hfill  !\r  
      !\hfill    -    !\hfill    -    !\hfill    -    !\hfill    -    !\hfill    -     !\hfill    -    !\hfill    25   !\hfill  260        !$yz$       !\hfill  !\rr  
\hfill {\bf 30}     !\hfill    1    !\hfill    4   !\hfill     11    !\hfill    25   !\hfill     50    !\hfill    91    !\hfill    154   !\hfill  246        !$xy+u^2,xz,$      !$H_{23}$  !$A_{23.1}$\r 
            !\hfill    -    !\hfill    -    !\hfill    -    !\hfill    1    !\hfill    9    !\hfill    46   !\hfill     175   !\hfill  550       !$x^2+z^2+u^2$       !\hfill  !\r  
            !\hfill    -    !\hfill    -    !\hfill    -    !\hfill    -    !\hfill    -     !\hfill    -    !\hfill    1   !\hfill  14        !$y^2,yu+zu$       !\hfill  !\rr  
\hfill {\bf 31}    !\hfill    1    !\hfill    4   !\hfill     11   !\hfill     25   !\hfill     50    !\hfill    91    !\hfill    154   !\hfill 246         !${\scriptstyle x^2-y^2,y^2-z^2,}$       !$H_{24}$ !$A_{23.1}$\r  
            !\hfill    -    !\hfill    -    !\hfill    -    !\hfill    2    !\hfill    16    !\hfill    77   !\hfill     282   !\hfill 864         !${\scriptstyle z^2-u^2,xz+yu,}$       !\hfill  !\r  
           !\hfill    -    !\hfill    -    !\hfill    -    !\hfill    -    !\hfill    -     !\hfill    -    !\hfill    4   !\hfill 48         !${\scriptstyle -x^2+xy-yz+xu}$       !\hfill !\rr   
\hfill {\bf 32}   !\hfill    1    !\hfill    4    !\hfill    11   !\hfill     25   !\hfill     51    !\hfill    97   !\hfill     176   !\hfill  309        !$x^2+z^2,xz,$      !$H_{23}$  !$A_{23.2}$\r  
            !\hfill    -    !\hfill    -    !\hfill    -    !\hfill    2   !\hfill     15    !\hfill    68   !\hfill     238   !\hfill 708         !$y^2,yu+zu,$       !\hfill !\r   
            !\hfill    -    !\hfill    -    !\hfill    -    !\hfill    -    !\hfill    -     !\hfill    -    !\hfill    4   !\hfill  44        !$u^2$       !\hfill !\rr   
\hfill {\bf 33}   !\hfill    1    !\hfill    4   !\hfill     11   !\hfill     25   !\hfill     51    !\hfill    97   !\hfill     176   !\hfill 309         !$x^2+xy,y^2+yz,$       !$H_{24}$   !$A_{23.2}$\r 
           !\hfill    -    !\hfill    -    !\hfill    -    !\hfill    3    !\hfill    22    !\hfill    99    !\hfill    345   !\hfill  1024        !$y^2+xu,z^2+xu,$       !\hfill  !\r  
           !\hfill    -    !\hfill    -    !\hfill    -    !\hfill    -    !\hfill    -     !\hfill    -    !\hfill    9   !\hfill 96         !$zu+u^2$       !\hfill  !\rr  
\hfill {\bf 34}    !\hfill    1    !\hfill    4   !\hfill     11   !\hfill     26    !\hfill    55    !\hfill    106   !\hfill     190   !\hfill 322         !${\scriptstyle x^2+xy+yu+u^2,y^2,xz,}$       !$H_{27}$   !$A_{27.1}$\r 
              !\hfill    -    !\hfill    -    !\hfill    -    !\hfill    -    !\hfill    -     !\hfill    5    !\hfill    38   !\hfill   172       !${\scriptstyle x^2+z^2+u^2,yu+zu}$       !\hfill !\rr   
\hfill {\bf 35}    !\hfill  1    !\hfill    4    !\hfill    11    !\hfill    26    !\hfill    55   !\hfill     108   !\hfill     201   !\hfill 360         !${\scriptstyle x^2+z^2+u^2,y^2,xz,}$   !$H_{27}$      !$A_{27.2}$\r 
    \hfill  !\hfill  -    !\hfill    -    !\hfill    -    !\hfill    -    !\hfill    2    !\hfill    16    !\hfill 76      !\hfill  278        !${\scriptstyle xy+yz+yu,yu+zu}$   !    !\rr   
\hfill {\bf 36}    ! \hfill 1  !  \hfill 4    !    \hfill 11   !\hfill   26   !\hfill  56   !\hfill 114     !\hfill 223   ! \hfill 424     !$x^2+y^2,z^2,u^2,$   !$H_{27}$    !$A_{27.3}$  \r 
\hfill  ! \hfill -  !  \hfill -    !    \hfill -   !\hfill   1   !\hfill  8   !\hfill 38     !\hfill 140   ! \hfill 441     !$yz-yu,$     ! !   \r 
\hfill  ! \hfill -  !  \hfill -    !    \hfill -   !\hfill   -   !\hfill  -   !\hfill -     !\hfill 1   !\hfill 12 !$xz+zu$ ! ! \rr 
\hfill {\bf 37} ! \hfill 1  !  \hfill 4    !    \hfill 11   !\hfill   26   !\hfill  56   !\hfill 114     !\hfill 223   ! \hfill 425     !$x^2,y^2,xy-zu,$   !$H_{27}$      !$A_{27.4}$   \r 
\hfill  ! \hfill -  !  \hfill -    !    \hfill -   !\hfill   1   !\hfill  8   !\hfill 38     !\hfill 141   ! \hfill 448     !$ yz-xu,$      ! !  \r 
\hfill  ! \hfill -  !  \hfill -    !    \hfill -   !\hfill   -   !\hfill  -   !\hfill -     !\hfill 1   ! \hfill 12    !$(x-y)(z-u)$     !  !   

\caption{\sl } 
$$

%
\font\bold=cmbx10 at 14pt

\newdimen\tempdim                
\newdimen\othick \othick=.4pt    
\newdimen\ithick \ithick=.4pt    
\newdimen\spacing \spacing=9pt   
\newdimen\abovehr \abovehr=6pt   
\newdimen\belowhr \belowhr=8pt   
\newdimen\nexttovr \nexttovr=8pt 

\def\r{\hfil&\omit\vrsp\vrule width\othick\cr&}   
\def\rr{\hfil\down{\abovehr}&\omit\vrsp\vrule width\othick\cr
 \noalign{\hrule height\ithick}\up{\belowhr}&}     
\def\up#1{\tempdim=#1\advance\tempdim by1ex
  \vrule height\tempdim width0pt depth0pt} 
\def\down#1{\vrule height0pt depth#1 width0pt} 
\def\large#1#2{\setbox0=\vtop{\hsize#1 \lineskiplimit=0pt \lineskip=1pt
  \baselineskip\spacing \advance\baselineskip by 3pt \noindent
  #2}\tempdim=\dp0\advance\tempdim by\abovehr\box0\down{\tempdim}}
\def\vrsp{\hskip\nexttovr\relax} 
\def\toprule#1{\def\startrule{\hrule height#1\relax}} 
\toprule{\othick}     
\def\nstrut{\vrule height\spacing depth3.5pt width0pt}
\def\exclaim{\char`\!}   
\def\preamble#1{\def\startup{#1}}  
\preamble{&##}  
{\catcode`\!=\active
\gdef!{\hfil\vrule width0pt\vrsp\vrule width\ithick\relax\vrsp&}}

\def\table #1{\vbox\bgroup \setbox0=\hbox{#1}
 \vbox\bgroup\offinterlineskip \catcode`\!=\active
\halign\bgroup##\vrule width\othick\vrsp&\span\startup\nstrut\cr
\noalign{\medskip}
\noalign{\startrule}\up{\belowhr}&}

\def\caption #1{\down{\abovehr}&\omit\vrsp\vrule width\othick\cr
\noalign{\hrule height\othick}\egroup\egroup \setbox1=\lastbox
\tempdim=\wd1 \hbox to\tempdim{\hfil \box0 \hfil} \box1 \smallskip
\hbox to\tempdim{\advance\tempdim by-20pt\hfil\vbox{\hsize\tempdim
\noindent #1}\hfil}\egroup}
$$\table{Table 5. Cases 38--57 }
Case !\multispan{8} The graded Betti Numbers! ideal gen:s (ex.)  ! $R(t)$ ! $R^{\exclaim}(t)$     \rr
\hfill{\bf 38}    !\hfill    1    !\hfill    4    !\hfill    11   !\hfill     26   !\hfill     57   !\hfill     120   !\hfill     247   !\hfill  502        !$x^2,y^2,z^2,zu,u^2$      !$H_{25}$   !$A_{25}$\rr  
\hfill {\bf 39}    !\hfill    1    !\hfill    4    !\hfill    11    !\hfill    26   !\hfill     57    !\hfill    120   !\hfill     247   !\hfill 502          !$x^2+yz+u^2,$      !$H_{26}$   !$A_{25}$\r  
      !\hfill    -    !\hfill    -    !\hfill    -    !\hfill    1    !\hfill    7    !\hfill    31   !\hfill     109     !\hfill  334       !$xz+z^2+yu,$       !\hfill  !\r  
       !\hfill    -    !\hfill    -    !\hfill    -    !\hfill    -    !\hfill    -     !\hfill    -    !\hfill    1   !\hfill 10         !$xy,xu,zu$       ! !\rr  
\hfill {\bf 40}    !\hfill    1    !\hfill    4    !\hfill    11    !\hfill    26   !\hfill     57   !\hfill     120    !\hfill    247   !\hfill 502         !$x^2-xu,xu-y^2,$      !$H_{27}$ !$A_{25}$\r   
     !\hfill   -   !\hfill    -    !\hfill    -    !\hfill    2    !\hfill     14    !\hfill    62   !\hfill     218     !\hfill  668       !$y^2-z^2,z^2-u^2,$       !\hfill !\r   
    !\hfill    -    !\hfill    -    !\hfill    -    !\hfill    -    !\hfill    -     !\hfill    -     !\hfill    4   !\hfill 40         !$xz+yu$      !\hfill !\rr   
\hfill {\bf 41}   !\hfill    1    !\hfill    4   !\hfill    11   !\hfill     27   !\hfill     62   !\hfill     137   !\hfill     295   !\hfill 624         !$xy,y^2,z^2,zu,u^2$      !$H_{29}$  !$A_{29}$\rr  
\hfill {\bf 42} ! \hfill 1 !\hfill 4     !\hfill 11   !\hfill 27   !\hfill 62   !\hfill 137     !\hfill 296  !\hfill 632   !$x^2+xy,zu,y^2,$   !$H_{29}$     !$A_{29.1}$  \r
\hfill  !\hfill -  !\hfill -    !\hfill -  !\hfill -   !\hfill -   !\hfill 1    !\hfill 8  !\hfill 30   !$xu,xz+yu$!     !  \rr
\hfill {\bf 43} ! \hfill 1 !\hfill 4     !\hfill 11   !\hfill 27   !\hfill 63   !\hfill 144     !\hfill 326  !\hfill 735  !$x^2,y^2,yz,zu,u^2$  !$H_{28}$      !$A_{28}$   \rr
\hfill {\bf 44}    !\hfill    1    !\hfill    4   !\hfill     11   !\hfill     27    !\hfill    63    !\hfill    144   !\hfill     326   !\hfill 735         !$xz,yz,y^2,$      !$H_{29}$  !$A_{28}$\r   
           !\hfill    -    !\hfill    -    !\hfill    -    !\hfill    1    !\hfill    7    !\hfill    31   !\hfill     111   !\hfill  352        !$yu+zu,$      !\hfill  !\r  
            !\hfill    -    !\hfill    -    !\hfill    -    !\hfill    -    !\hfill    -     !\hfill    -    !\hfill    1   !\hfill 10         !$z^2+u^2$       !   !\rr  
\hfill {\bf 45} ! \hfill 1 !\hfill 4     !\hfill 11   !\hfill 28   !\hfill 68   !\hfill 162     !\hfill 382  !\hfill 896   !${\scriptstyle xy+yz,xy+z^2+yu,}$!$H_{31}$      !$A_{31}$ \r
\hfill  ! \hfill !\hfill     !\hfill   !\hfill  !\hfill  !\hfill      !\hfill  !\hfill    !${\scriptstyle yu+zu,y^2,xz}$!     !  \rr
\hfill {\bf 46} ! \hfill 1 !\hfill 4     !\hfill 11   !\hfill 28   !\hfill 69   !\hfill 168     !\hfill 407  !\hfill  984       !$x^2,xy,yz,zu,u^2$ !$H_{30}$  !$A_{30}$  \rr
\hfill {\bf 47} ! \hfill 1 !\hfill 4     !\hfill 11   !\hfill 28   !\hfill 69   !\hfill 168     !\hfill 407  !\hfill 984  !$x^2+xy,y^2,$      !$H_{31}$  !$A_{30}$   \r
\hfill    ! \hfill -  !  \hfill -    !    \hfill -    !\hfill    1   !\hfill   6   !\hfill 25   !\hfill 88   ! \hfill 282 !$xu,xz+yu,$        !  ! \r 
\hfill    ! \hfill -  !  \hfill -    !    \hfill -    !\hfill   -    !\hfill   -   !\hfill -    !\hfill   1  !\hfill 8    !$-x^2+xz-yz$       !  !  \rr
\hfill {\bf 48}    !\hfill    1    !\hfill    4   !\hfill     11   !\hfill     28   !\hfill     69   !\hfill     169   !\hfill     413   !\hfill 1009  !$xy,z^2+yu,$      !$H_{31}$  !$A_{31.1}$\r 
        !\hfill    -    !\hfill    -    !\hfill    -    !\hfill    1    !\hfill    7    !\hfill    31   !\hfill     113   !\hfill 370         !$yu+zu,y^2,$       !\hfill   !\r 
         !\hfill    -    !\hfill    -    !\hfill    -    !\hfill    -    !\hfill    -     !\hfill    -    !\hfill    1     !\hfill 10   !$xz$   !\hfill    !\rr
\hfill  {\bf 49}   !\hfill    1    !\hfill    4    !\hfill    11   !\hfill     29   !\hfill     75   !\hfill     193   !\hfill     496     !\hfill 1274      !$xz,y^2,z^2,yu,zu$ !$H_{33}$ !$A_{33}$\rr   
\hfill 50    !\hfill    1    !\hfill    4    !\hfill    11   !\hfill     29   !\hfill     76    !\hfill    199   !\hfill     521   !\hfill     1364     !$x^2,xy,xz,y^2,z^2$        !$H_{32}$ !$A_{32}$\rr   
\hfill {\bf 51 }   !\hfill    1    !\hfill    4    !\hfill    11    !\hfill    29    !\hfill    76   !\hfill    199   !\hfill     521   !\hfill     1364     !$xy,xz,yz+xu,$       !$H_{33}$!$A_{32}$\r   
      !\hfill    -    !\hfill    -    !\hfill    -    !\hfill    1    !\hfill    6    !\hfill    25    !\hfill    90    !\hfill    300     !$z^2,zu$       !\hfill !\r   
       !\hfill    -    !\hfill    -    !\hfill    -    !\hfill    -    !\hfill    -     !\hfill    -     !\hfill    1     !\hfill   8  !\hfill              !\hfill !\rr   
\hfill 52    !\hfill    1    !\hfill    4    !\hfill    11   !\hfill     30   !\hfill     82   !\hfill     224    !\hfill    612     !\hfill 1672        !$x^2,xy,xz,y^2,yz$       !$H_{34}$ !$A_{34}$\rr    
\hfill {\bf 53}    !\hfill    1    !\hfill    4    !\hfill    11    !\hfill    31   !\hfill     88   !\hfill     249   !\hfill     705     !\hfill 1996        !${\scriptstyle y^2-u^2,xz,yz,z^2,zu}$       !$H_{35}$ !$A_{35}$\rr   
\hfill{\bf  54}    !\hfill    1    !\hfill    4    !\hfill    12    !\hfill    32   !\hfill     80   !\hfill     192    !\hfill    448    !\hfill    1024   !${\scriptstyle x^2,xz,y^2,z^2,yu+zu,u^2}$! $H_{13}$       !$A_{13}$ \rr  
\hfill {\bf  55}    !\hfill    1    !\hfill    4    !\hfill    12    !\hfill    32    !\hfill    80    !\hfill    192    !\hfill    449    !\hfill    1034     !$x^2+xy,xz+yu,$      ! $H_{13}$   !$A_{13.1}$\r  
      !\hfill    -    !\hfill    -    !\hfill    -    !\hfill    -    !\hfill    -     !\hfill    1    !\hfill    10     !\hfill    57     !$xu,y^2,z^2,zu+u^2$          !\hfill       !\hfill  \rr  
\hfill {\bf 56} ! \hfill 1   !\hfill 4     !\hfill 12       !\hfill 32   !\hfill 80   !\hfill 193     !\hfill 457  !\hfill 1072   !$x^2+xz+u^2,xy,$ !  $H_{13}$     !$A_{13.2}$  \r
\hfill    ! \hfill -  !  \hfill -    !    \hfill -   !\hfill   -   !\hfill   1   !\hfill 9     !\hfill 48   ! \hfill 199    !$xu,x^2-y^2,z^2,zu$     !  ! \rr 
 \hfill{\bf 57 }! \hfill 1   !\hfill 4     !\hfill 12       !\hfill 32   !\hfill 81   !\hfill 200     !\hfill 488  !\hfill 1184   !${\scriptstyle x^2+yz+u^2,xu,}$!$H_{13}$      ! $A_{13.3}$  \r
\hfill    ! \hfill -  !  \hfill -    !    \hfill -   !\hfill   1   !\hfill   8   !\hfill 40     !\hfill 160   ! \hfill 562    !${\scriptstyle x^2+xy,xz+yu,}$     !  !  \r 
\hfill    ! \hfill -  !  \hfill -    !    \hfill -    !\hfill   -    !\hfill   -   !\hfill -      !\hfill   1   !\hfill 12 !${\scriptstyle zu+u^2,y^2+z^2}$ ! !
\caption{\sl } 
$$
%
\font\bold=cmbx10 at 14pt

\newdimen\tempdim                
\newdimen\othick \othick=.4pt    
\newdimen\ithick \ithick=.4pt    
\newdimen\spacing \spacing=9pt   
\newdimen\abovehr \abovehr=6pt   
\newdimen\belowhr \belowhr=8pt   
\newdimen\nexttovr \nexttovr=8pt 

\def\r{\hfil&\omit\vrsp\vrule width\othick\cr&}   
\def\rr{\hfil\down{\abovehr}&\omit\vrsp\vrule width\othick\cr
 \noalign{\hrule height\ithick}\up{\belowhr}&}     
\def\up#1{\tempdim=#1\advance\tempdim by1ex
  \vrule height\tempdim width0pt depth0pt} 
\def\down#1{\vrule height0pt depth#1 width0pt} 
\def\large#1#2{\setbox0=\vtop{\hsize#1 \lineskiplimit=0pt \lineskip=1pt
  \baselineskip\spacing \advance\baselineskip by 3pt \noindent
  #2}\tempdim=\dp0\advance\tempdim by\abovehr\box0\down{\tempdim}}
\def\vrsp{\hskip\nexttovr\relax} 
\def\toprule#1{\def\startrule{\hrule height#1\relax}} 
\toprule{\othick}     
\def\nstrut{\vrule height\spacing depth3.5pt width0pt}
\def\exclaim{\char`\!}   
\def\preamble#1{\def\startup{#1}}  
\preamble{&##}  
{\catcode`\!=\active
\gdef!{\hfil\vrule width0pt\vrsp\vrule width\ithick\relax\vrsp&}}

\def\table #1{\vbox\bgroup \setbox0=\hbox{#1}
 \vbox\bgroup\offinterlineskip \catcode`\!=\active
\halign\bgroup##\vrule width\othick\vrsp&\span\startup\nstrut\cr
\noalign{\medskip}
\noalign{\startrule}\up{\belowhr}&}

\def\caption #1{\down{\abovehr}&\omit\vrsp\vrule width\othick\cr
\noalign{\hrule height\othick}\egroup\egroup \setbox1=\lastbox
\tempdim=\wd1 \hbox to\tempdim{\hfil \box0 \hfil} \box1 \smallskip
\hbox to\tempdim{\advance\tempdim by-20pt\hfil\vbox{\hsize\tempdim
\noindent #1}\hfil}\egroup}
$$\table{Table 6. Cases 58--74 }
Case !\multispan{8} The graded Betti Numbers! ideal gen:s,ex. ! $R(t)$ ! $R^{\exclaim}(t)$     \rr
\hfill {\bf 58 }! \hfill 1  !  \hfill 4    !    \hfill 12   !\hfill  33   !\hfill  87   !\hfill 225     !\hfill 576   !\hfill 1467     !${\scriptstyle x^2+xy,x^2+zu,y^2,}$! $H_{15}$     !$A_{15}$   \r 
\hfill ! \hfill  !  \hfill    !    \hfill    !\hfill     !\hfill     !\hfill      !\hfill   !\hfill      !${\scriptstyle z^2,xz+yu,xu}$!  !  \rr 
\hfill {\bf 59} ! \hfill 1  !  \hfill 4    !    \hfill 12   !\hfill  33   !\hfill  87   !\hfill 225     !\hfill 576   ! \hfill 1468     !${\scriptstyle x^2-y^2,xy,xu,}$   ! $H_{15}$   ! $A_{15.1}$ \r 
\hfill    ! \hfill -  !  \hfill -    !    \hfill -    !\hfill    -   !\hfill   -   !\hfill -     !\hfill 1   ! \hfill 8     !${\scriptstyle z^2,zu,xz+yu}$     !  !  \rr 
\hfill {\bf 60} ! \hfill 1  !  \hfill 4    !    \hfill 12   !\hfill   33   !\hfill  87   !\hfill 225     !\hfill 577   ! \hfill 1474     !${\scriptstyle x^2+yz+u^2,xz+yu,}$   ! $H_{15} $   !$A_{15.2} $  \r 
\hfill    ! \hfill -  !  \hfill -    !    \hfill -    !\hfill    -   !\hfill   -   !\hfill 1     !\hfill 7   ! \hfill 33     !${\scriptstyle zu,x^2+xy,z^2,xu}$   !  !  \rr 
\hfill {\bf 61} ! \hfill 1  !  \hfill 4    !    \hfill 12   !\hfill  33   !\hfill  88   !\hfill 232     !\hfill 609   !\hfill 1596     !${\scriptstyle x^2-y^2,xy,z^2,}$! $H_{14}$ !$A_{14}$   \r 
\hfill  ! \hfill  !  \hfill     !    \hfill    !\hfill    !\hfill     !\hfill      !\hfill    !\hfill      !${\scriptstyle xu,zu,u^2}$!  !  \rr 
\hfill {\bf 62} ! \hfill 1  !  \hfill 4    !    \hfill 12   !\hfill   33   !\hfill  88   !\hfill 232     !\hfill 609   ! \hfill 1596     !$x^2-y^2,xy,$   ! $H_{15} $   !$A_{14}$  \r 
\hfill    ! \hfill -  !  \hfill -    !    \hfill -    !\hfill    1   !\hfill   7   !\hfill 33     !\hfill 129   ! \hfill 455     !$xu,yz+yu,$     !  !  \r 
\hfill    ! \hfill -  !  \hfill -    !    \hfill -    !\hfill   -    !\hfill   -   !\hfill -      !\hfill   1   !\hfill 10 !$z^2,zu$ ! !  \rr
\hfill {\bf 63} ! \hfill 1  !  \hfill 4    !    \hfill 12   !\hfill  34   !\hfill  94   !\hfill 258     !\hfill 706   !\hfill 1930     !${\scriptstyle x^2,xy,xu,y^2,z^2,zu}$     !$H_{17}$    ! $A_{17}$ \rr 
\hfill {\bf 64} ! \hfill 1  !  \hfill 4    !    \hfill 12   !\hfill  34   !\hfill  95   !\hfill 265     !\hfill 739   !\hfill 2061     !${\scriptstyle x^2-y^2,xy, }$     ! $H_{16}$  !$A_{16}$ \r 
\hfill ! \hfill   !  \hfill      !    \hfill     !\hfill     !\hfill      !\hfill       !\hfill     !\hfill       !${\scriptstyle  z^2,xu,yu,zu}$     !   ! \rr 
\hfill {\bf 65} ! \hfill 1  !  \hfill 4    !    \hfill 12   !\hfill  34   !\hfill  95   !\hfill 265     !\hfill 739   !\hfill 2061     !$x^2,xy,xz,$     !$H_{17}$   !$A_{16}$  \r 
\hfill    ! \hfill -  !  \hfill -    !    \hfill -    !\hfill   1   !\hfill   7   !\hfill 33     !\hfill 131   ! \hfill 475     !$y^2,yu+z^2,$     !  !  \r 
\hfill    ! \hfill -  !  \hfill -    !    \hfill -    !\hfill   -    !\hfill   -   !\hfill -      !\hfill   1   !\hfill 10 !$yu+zu$ ! !  \rr
\hfill {\bf 66} ! \hfill 1  !  \hfill 4    !    \hfill 12   !\hfill  35   !\hfill  101   !\hfill 291     !\hfill 838   !\hfill 2413     !${\scriptstyle xz,y^2,yu,z^2,zu,u^2}$     ! $H_{18}$  !$A_{18}$ \rr 
\hfill {\bf 67} ! \hfill 1  !  \hfill 4    !    \hfill 12   !\hfill  36   !\hfill  107   !\hfill 318     !\hfill 945   ! \hfill 2808     !${\scriptstyle xy,xz,y^2,yu,z^2,zu}$     ! $H_{20}$  !$A_{20}$ \rr 
\hfill 68 ! \hfill 1  !  \hfill 4    !    \hfill 12   !\hfill  36   !\hfill  108   !\hfill 324     !\hfill 972   !\hfill 2916      !${\scriptstyle x^2,xy,xz,y^2,yz,z^2}$     ! $H_{19}$  !$A_{19}$  \rr 
\hfill {\bf 69 }! \hfill 1  !  \hfill 4    !    \hfill 12   !\hfill  36   !\hfill  108   !\hfill 324     !\hfill 972   !\hfill 2916      !$x^2,xz,xu,$     !$H_{20}$   !$A_{19}$ \r 
\hfill    ! \hfill -  !  \hfill -    !    \hfill -    !\hfill   1   !\hfill   6    !\hfill 27     !\hfill 108   ! \hfill 405     !$xy-zu,$     !   ! \r 
\hfill    ! \hfill -  !  \hfill -    !    \hfill -    !\hfill   -    !\hfill   -   !\hfill -      !\hfill   1   !\hfill 8 !$yz,z^2$ !  ! \rr
\hfill {\bf 70} ! \hfill 1  !  \hfill 4    !    \hfill 12   !\hfill  37   !\hfill  114   !\hfill 351     !\hfill 1081   !\hfill 3329      !${\scriptstyle x^2,xy,xz,xu,y^2,yz}$       !$H_{21}$ !$A_{21}$ \rr 
\hfill {\bf 71} ! \hfill 1  !  \hfill 4    !    \hfill 13   !\hfill  40   !\hfill  121   !\hfill 364     !\hfill 1093   !\hfill 3280     !${\scriptstyle x^2,y^2,z^2,u^2,}$       !$H_{7}$ !$A_{7}$ \r
\hfill          ! \hfill    !  \hfill      !    \hfill      !\hfill       !\hfill        !\hfill         !\hfill        !\hfill          !${\scriptstyle xy,zu,yz+xu}$       ! ! \rr   
\hfill {\bf 72 }! \hfill 1  !  \hfill 4    !    \hfill 13   !\hfill  41   !\hfill  128   !\hfill 399     !\hfill 1243   !\hfill 3872      !${\scriptstyle x^2-y^2,xy,yz,}$       ! ! \r 
\hfill   ! \hfill    !  \hfill     !    \hfill   !\hfill     !\hfill   !\hfill      !\hfill     !      !${\scriptstyle zu,z^2,xz+yu,xu}$       !$H_{9}$ !$A_{9}$ \rr 
\hfill {\bf 73} ! \hfill 1  !  \hfill 4    !    \hfill 13   !\hfill  41   !\hfill  129   !\hfill 406     !\hfill 1278   ! \hfill 4023     !${\scriptstyle x^2,y^2,z^2,u^2,}$       !$H_{8}$ !$A_{8}$ \r 
\hfill    ! \hfill    !  \hfill      !    \hfill      !\hfill       !\hfill        !\hfill         !\hfill        !      !${\scriptstyle zu,yu,xu}$       ! ! \rr 
\hfill {\bf 74}! \hfill 1  !  \hfill 4    !    \hfill 13   !\hfill  41   !\hfill  129   !\hfill 406     !\hfill 1278   !\hfill 4023  !$x^2,xy+z^2,$       !$H_9$ !$A_{8}$ \r 
\hfill    ! \hfill -  !  \hfill -    !    \hfill -    !\hfill    1   !\hfill   7   !\hfill 35     !\hfill 151   ! \hfill 604  !$yz,xu,yu,$  ! ! \r 
\hfill    ! \hfill -  !  \hfill -    !    \hfill -    !\hfill   -    !\hfill   -   !\hfill -      !\hfill   1   !\hfill  10 ! $zu,u^2$ ! ! 
\caption{\sl } $$
%
\font\bold=cmbx10 at 14pt

\newdimen\tempdim                
\newdimen\othick \othick=.4pt    
\newdimen\ithick \ithick=.4pt    
\newdimen\spacing \spacing=9pt   
\newdimen\abovehr \abovehr=6pt   
\newdimen\belowhr \belowhr=8pt   
\newdimen\nexttovr \nexttovr=8pt 

\def\r{\hfil&\omit\vrsp\vrule width\othick\cr&}   
\def\rr{\hfil\down{\abovehr}&\omit\vrsp\vrule width\othick\cr
 \noalign{\hrule height\ithick}\up{\belowhr}&}     
\def\up#1{\tempdim=#1\advance\tempdim by1ex
  \vrule height\tempdim width0pt depth0pt} 
\def\down#1{\vrule height0pt depth#1 width0pt} 
\def\large#1#2{\setbox0=\vtop{\hsize#1 \lineskiplimit=0pt \lineskip=1pt
  \baselineskip\spacing \advance\baselineskip by 3pt \noindent
  #2}\tempdim=\dp0\advance\tempdim by\abovehr\box0\down{\tempdim}}
\def\vrsp{\hskip\nexttovr\relax} 
\def\toprule#1{\def\startrule{\hrule height#1\relax}} 
\toprule{\othick}     
\def\nstrut{\vrule height\spacing depth3.5pt width0pt}
\def\exclaim{\char`\!}   
\def\preamble#1{\def\startup{#1}}  
\preamble{&##}  
{\catcode`\!=\active
\gdef!{\hfil\vrule width0pt\vrsp\vrule width\ithick\relax\vrsp&}}

\def\table #1{\vbox\bgroup \setbox0=\hbox{#1}
 \vbox\bgroup\offinterlineskip \catcode`\!=\active
\halign\bgroup##\vrule width\othick\vrsp&\span\startup\nstrut\cr
\noalign{\medskip}
\noalign{\startrule}\up{\belowhr}&}

\def\caption #1{\down{\abovehr}&\omit\vrsp\vrule width\othick\cr
\noalign{\hrule height\othick}\egroup\egroup \setbox1=\lastbox
\tempdim=\wd1 \hbox to\tempdim{\hfil \box0 \hfil} \box1 \smallskip
\hbox to\tempdim{\advance\tempdim by-20pt\hfil\vbox{\hsize\tempdim
\noindent #1}\hfil}\egroup}
$$\table{Table 7. Cases 75--83 }
Case !\multispan{8} The graded Betti Numbers! ideal gen:s ! $R(t)$ ! $R^{\exclaim}(t)$     \rr
\hfill {\bf 75} ! \hfill 1  !  \hfill 4    !    \hfill 13   !\hfill  42   !\hfill  135   !\hfill 434     !\hfill 1395   !\hfill 4484     !${\scriptstyle x^2,xy,xz,xu, }$       !$H_{10}$ !$A_{10}$ \r 
\hfill    ! \hfill    !  \hfill      !    \hfill      !\hfill       !\hfill        !\hfill         !\hfill    !      !${\scriptstyle  y^2,yz,u^2}$       ! ! \rr 
\hfill {\bf 76} ! \hfill 1  !  \hfill 4    !    \hfill 13   !\hfill  43   !\hfill  142   !\hfill 469     !\hfill 1549   !\hfill 5116      !${\scriptstyle x^2,xy,xz,xu,}$     !$H_{11}$   !$A_{11}$ \r 
\hfill    ! \hfill    !  \hfill      !    \hfill      !\hfill       !\hfill        !\hfill         !\hfill    !      !${\scriptstyle  z^2,zu,yu}$     !   ! \rr 
\hfill {\bf 77} ! \hfill 1  !  \hfill 4    !    \hfill 13   !\hfill  44   !\hfill  148   !\hfill 498     !\hfill 1676   !\hfill 5640      !${\scriptstyle x^2,xy,xz,xu,}$       !$H_{12}$ !$A_{12}$ \r 
\hfill    ! \hfill    !  \hfill      !    \hfill      !\hfill       !\hfill        !\hfill         !\hfill        !      !${\scriptstyle y^2,yz,yu}$       ! ! \rr 
\hfill {\bf 78} ! \hfill 1  !  \hfill 4    !    \hfill 14   !\hfill  48   !\hfill  164   !\hfill 560     !\hfill 1912   !\hfill 6528      !${\scriptstyle x^2,xy,y^2,z^2,zu, }$!$H_{4}$ !$A_{4}$ \r 
\hfill    ! \hfill    !  \hfill      !    \hfill      !\hfill       !\hfill        !\hfill         !\hfill        !      !${\scriptstyle u^2,xz+yu,yz-xu }$! ! \rr 
\hfill {\bf 79} ! \hfill 1  !  \hfill 4    !    \hfill 14   !\hfill  49   !\hfill  171   !\hfill 597     !\hfill 2084   !\hfill 7275      !${\scriptstyle x^2,xy,xz,xu, }$       !$H_{5}$  !$A_{5}$ \r 
\hfill    ! \hfill    !  \hfill      !    \hfill      !\hfill       !\hfill        !\hfill         !\hfill        !      !${\scriptstyle y^2,yu,z^2,zu }$       ! ! \rr 
\hfill {\bf 80} ! \hfill 1  !  \hfill 4    !    \hfill 14   !\hfill  50   !\hfill  178   !\hfill 634     !\hfill 2258   ! \hfill 8042     !${\scriptstyle x^2,xy,xz,y^2,}$       !$H_{6}$ !$A_{6}$ \r
\hfill    ! \hfill    !  \hfill      !    \hfill      !\hfill       !\hfill        !\hfill         !\hfill        !      !${\scriptstyle yz,yu,z^2,zu }$        ! ! \rr 
\hfill {\bf 81} ! \hfill 1  !  \hfill 4    !    \hfill 15  !\hfill   56   !\hfill  209   !\hfill 780     !\hfill 2911   ! \hfill 10864     !${\scriptstyle x^2,y^2,z^2,u^2,xy,}$     !$H_{2}$   !$A_{2}$ \r 
\hfill   ! \hfill    !  \hfill      !    \hfill    !\hfill       !\hfill      !\hfill      !\hfill    !      !${\scriptstyle xz,yz-xu,yu,zu}$     !   ! \rr 
\hfill {\bf 82} ! \hfill 1  !  \hfill 4    !    \hfill 15  !\hfill   57   !\hfill  216   !\hfill  819    !\hfill 3105   ! 11772     !${\scriptstyle x^2,xy,xz,xu,y^2,}$     !$H_3$   !$A_3$ \r 
\hfill   ! \hfill    !  \hfill      !    \hfill    !\hfill       !\hfill      !\hfill      !\hfill    !      !${\scriptstyle zu,u^2,yz,yu}$     !   ! \rr 
\hfill {\bf 83} ! \hfill 1  !  \hfill 4    !    \hfill 16   !\hfill   64   !\hfill  256   !\hfill 1024     !\hfill 4096   ! \hfill 16384           !${\scriptstyle x^2,y^2,z^2,u^2,xy,}$ !$H_1$ !$A_1$\r
\hfill  ! \hfill  !  \hfill     !    \hfill    !\hfill      !\hfill    !\hfill      !\hfill    ! \hfill            !${\scriptstyle xz,xu,yz,yu,zu}$ ! !
\caption{\sl } 
$$
%
\font\bold=cmbx10 at 14pt

\newdimen\tempdim                
\newdimen\othick \othick=.4pt    
\newdimen\ithick \ithick=.4pt    
\newdimen\spacing \spacing=9pt   
\newdimen\abovehr \abovehr=6pt   
\newdimen\belowhr \belowhr=8pt   
\newdimen\nexttovr \nexttovr=8pt 

\def\r{\hfil&\omit\vrsp\vrule width\othick\cr&}   
\def\rr{\hfil\down{\abovehr}&\omit\vrsp\vrule width\othick\cr
 \noalign{\hrule height\ithick}\up{\belowhr}&}     
\def\up#1{\tempdim=#1\advance\tempdim by1ex
  \vrule height\tempdim width0pt depth0pt} 
\def\down#1{\vrule height0pt depth#1 width0pt} 
\def\large#1#2{\setbox0=\vtop{\hsize#1 \lineskiplimit=0pt \lineskip=1pt
  \baselineskip\spacing \advance\baselineskip by 3pt \noindent
  #2}\tempdim=\dp0\advance\tempdim by\abovehr\box0\down{\tempdim}}
\def\vrsp{\hskip\nexttovr\relax} 
\def\toprule#1{\def\startrule{\hrule height#1\relax}} 
\toprule{\othick}     
\def\nstrut{\vrule height\spacing depth3.5pt width0pt}
\def\exclaim{\char`\!}   
\def\preamble#1{\def\startup{#1}}  
\preamble{&##}  
{\catcode`\!=\active
\gdef!{\hfil\vrule width0pt\vrsp\vrule width\ithick\relax\vrsp&}}

\def\table #1{\vbox\bgroup \setbox0=\hbox{#1}
 \vbox\bgroup\offinterlineskip \catcode`\!=\active
\halign\bgroup##\vrule width\othick\vrsp&\span\startup\nstrut\cr
\noalign{\medskip}
\noalign{\startrule}\up{\belowhr}&}

\def\caption #1{\down{\abovehr}&\omit\vrsp\vrule width\othick\cr
\noalign{\hrule height\othick}\egroup\egroup \setbox1=\lastbox
\tempdim=\wd1 \hbox to\tempdim{\hfil \box0 \hfil} \box1 \smallskip
\hbox to\tempdim{\advance\tempdim by-20pt\hfil\vbox{\hsize\tempdim
\noindent #1}\hfil}\egroup}
$$\table{Table 8. The 16 series $R^!(t)$ that do not correspond to Koszul algebras. }
 Rings ! $R^{\exclaim}(t)$ ! Numerator of $R^{\exclaim}(t)$ ! Denominator of $R^{\exclaim}(t)$   \rr
\hfill $17$ ! $A_{43.1}$    ! $(1+t)^4(1+t^3)$    !$(1-t^2)^4$ \rr
\hfill $22$ ! $A_{46.1}$    ! $(1+t)^2(1-t+t^2)$!$(1-t)^2(1-t-t^2-t^4-t^5)$  \rr

\hfill 29     !  $A_{22.1}$         !   $(1+t)^4$ !$(1-t^2)^5$       \rr

\hfill 30,31  !  $A_{23.1}$      !      $(1+t)^4(1+t^3)$ !$(1-t^2)^5$  \rr

\hfill 32,33  !  $A_{23.2}$      !     $(1+t)^3$ !$(1-t^2)^3(1-t-t^2)$   \rr

\hfill 34    !   $A_{27.1}$      !     $(1+t)^4(1+t^3)^2$ !$(1-t^2)^5 (1-t^4)$  \rr

\hfill 35    !   $A_{27.2}$      !     $(1+t)^3(1+t^3)$ ! $(1-t^2)^3 (1-t-t^2)$  \rr

\hfill 36    !   $A_{27.3}$      !    $(1+t)^2$ !$(1-t^2) (1-t-t^2)^2$     \rr

\hfill 37    !   $A_{27.4}$      !   $(1+t)^3(1+t^3)$ !$(1-t^2)^3(1-t-t^2-t^4-t^5)$ \rr

\hfill 42    !   $A_{29.1}$      !  $(1+t)^2(1+t^3)$! $(1-t^2)^2(1-2\,t-t^4)$ \rr

\hfill 48    !   $A_{31.1}$  !$(1+t)^2$ !$1-2\,t-2\,t^2+2\,t^3+t^4-t^5$   \rr

\hfill 55    !   $A_{13.1}$   !$(1-t+t^2)^2$ !$(1-t)^3(1-3\,t+3\,t^2-3\,t^3)$  \rr

\hfill 56    !   $A_{13.2}$    !$(1-t+t^2)$ !$(1-t)^2(1-3\,t+2\,t^2-t^3)$ \rr

\hfill 57    !   $A_{13.3}$    !$1$ !$(1-t)^2(1-2\,t-t^2)$\rr

\hfill 59    !   $A_{15.1}$   !$(1+t^3)$!$(1-t)^2(1-2\,t-t^2-2\,t^4-t^5)$ \rr

\hfill 60    !   $A_{15.2}$  !$1$ !$(1-t)^2(1-2\,t-t^2-t^3)$  
\caption{\sl For e.g. the ring $34$ the nilpotency degree of $\eta$ is $5$.} 
$$
 \end